\documentclass[pdflatex,sn-mathphys-num]{sn-jnl}

\usepackage{amsfonts,amssymb,mathtools,graphicx}
\usepackage{booktabs,array,microtype}

\usepackage{xcolor} % 如果已经加载，省略这一行

\usepackage{amsthm}
\usepackage{aliascnt}
\usepackage{tikz}
\usepackage{caption}
\theoremstyle{thmstyleone}
\newtheorem{theorem}{Theorem}[section]

\newaliascnt{lemma}{theorem}
\newtheorem{lemma}[lemma]{Lemma}
\aliascntresetthe{lemma}

\newaliascnt{proposition}{theorem}
\newtheorem{proposition}[proposition]{Proposition}
\aliascntresetthe{proposition}

\newaliascnt{corollary}{theorem}
\newtheorem{corollary}[corollary]{Corollary}
\aliascntresetthe{corollary}

\theoremstyle{thmstylethree}
\newaliascnt{remark}{theorem}
\newtheorem{remark}[remark]{Remark}
\aliascntresetthe{remark}

\usepackage{hyperref}
\usepackage[nameinlink,capitalize,noabbrev]{cleveref}

\crefname{theorem}{Theorem}{Theorems}
\Crefname{theorem}{Theorem}{Theorems}

\crefname{lemma}{Lemma}{Lemmas}
\Crefname{lemma}{Lemma}{Lemmas}

\crefname{proposition}{Proposition}{Propositions}
\Crefname{proposition}{Proposition}{Propositions}

\crefname{corollary}{Corollary}{Corollaries}
\Crefname{corollary}{Corollary}{Corollaries}

\crefname{remark}{Remark}{Remarks}
\Crefname{remark}{Remark}{Remarks}
\crefname{equation}{}{}
\Crefname{equation}{}{}

\allowdisplaybreaks[2]
\numberwithin{equation}{section}
\newcommand{\R}{\mathbb R}
\newcommand{\C}{\mathbb C}
\newcommand{\abs}[1]{\lvert#1\rvert}
\newcommand{\one}{\mathbf 1}
\newcommand{\ind}{\mathbf 1}
\newcommand{\dd}{\,\mathrm d}
\newcommand{\ii}{\mathrm i}
\newcommand{\wh}{\widehat}
\DeclareMathOperator{\Tr}{Tr}
\DeclareMathOperator{\Ran}{Ran}
\DeclareMathOperator{\Ker}{Ker}
\DeclareMathOperator{\rank}{rank}
\newcommand{\HS}{\mathrm{HS}}
\title[Improved Berezin--Li--Yau bounds]
{Improved Berezin--Li--Yau Bounds for Dirichlet Eigenvalues via
Fourier Coupling and Frequency Concentration}

\author[1]{\fnm{Yifan} \sur{Wang}}\email{wangyifan@lsec.cc.ac.cn}
\author*[2,3]{\fnm{Hehu} \sur{Xie}}\email{hhxie@lsec.cc.ac.cn}

\affil[1]{\orgdiv{School of Statistics and Mathematics},
\orgname{Central University of Finance and Economics},
\orgaddress{\city{Beijing}, \postcode{102206}, \country{China}}}

\affil*[2]{\orgdiv{SKLMS, NCMIS, Institute of Computational Mathematics},
\orgname{Academy of Mathematics and Systems Science, Chinese Academy of Sciences},
\orgaddress{\street{No.~55 Zhongguancun Donglu}, \city{Beijing},
\postcode{100190}, \country{China}}}

\affil[3]{\orgdiv{School of Mathematical Sciences},
\orgname{University of Chinese Academy of Sciences},
\orgaddress{\city{Beijing}, \postcode{100049}, \country{China}}}

\begin{document}
\abstract{Let $0<\lambda_1\leq\lambda_2\leq\cdots$ be the Dirichlet
eigenvalues of an open set $\Omega\subset\mathbb R^n$, $n\geq2$.  We
develop two Fourier-space refinements of the Berezin--Li--Yau argument.
For every finite-measure open set, coordinate commutator identities and
coupled scalar--vector projections lead to a one-dimensional
variational constant $C_n^{\mathrm{new}}$ and the individual bound
\[
\lambda_k\geq C_n^{\mathrm{new}}(2\pi)^2\omega_n^{-2/n}
|\Omega|^{-2/n}k^{2/n},
\qquad
\frac{n}{n+2}<c_n^{\mathrm{env}}<C_n^{\mathrm{new}}<1.
\]
Here $n/(n+2)$ is the classical Li--Yau individual coefficient and
$c_n^{\mathrm{env}}$ is the coefficient associated with the uncoupled
scalar envelope developed in our previous paper. 
Thus the displayed strict hierarchy separates the gain of the scalar
energy envelope from the further gain supplied by the commutator
coupling.  For bounded open sets, a frequency-concentration
operator, its first two traces, and the Ky Fan principle force Fourier
mass outside the Li--Yau filling ball.  A two-region bathtub argument
then gives
\[
\sum_{j=1}^k\lambda_j
\geq \mathcal C_k^{(2)}(\Omega)
\frac{n}{n+2}(2\pi)^2\omega_n^{-2/n}
|\Omega|^{-2/n}k^{1+2/n},
\qquad \mathcal C_k^{(2)}(\Omega)>1,
\]
where $\mathcal C_k^{(2)}(\Omega)$ is strictly decreasing in $k$ and
converges to one, and
\[
\mathcal C_k^{(2)}(\Omega)-1\geq
\frac{n+2}{2(n+1)^2}
\left(\frac{\omega_{n-1}}{\omega_n}\right)^2
\left(\frac34\right)^{n-1}k^{-2/n}.
\]
Thus the sharp leading term in the Berezin--Li--Yau sum is supplemented
by an explicit positive, volume-only correction of order $k$.  Apart
from boundedness in the concentration argument, no boundary regularity,
convexity, or connectedness is assumed.  These improvements do not
prove P\'{o}lya's conjecture: the individual coefficient is below one,
and the sum multiplier tends to one.}


\maketitle

\section{Introduction}\label{sec:introduction}

Let $\Omega\subset\R^n$, $n\geq2$, be an open set of finite positive
measure $V=|\Omega|$, and let
\[
0<\lambda_1\leq\lambda_2\leq\cdots
\]
be the Dirichlet eigenvalues, repeated according to multiplicity.  With
$\omega_n=|B_1(0)|$, set
\begin{equation}\label{eq:weyl-constant}
C_n=(2\pi)^2\omega_n^{-2/n}.
\end{equation}
Weyl's law states that
$\lambda_k\sim C_nV^{-2/n}k^{2/n}$ as $k\to\infty$; see, for example,
\cite{SafarovVassiliev1997}.  P\'olya's conjecture~\cite{Polya1961}
asks whether the corresponding termwise inequality
\begin{equation}\label{eq:polya}
\lambda_k\geq C_nV^{-2/n}k^{2/n}
\end{equation}
holds for every $k$.  This remains open for general domains.

The classical Berezin--Li--Yau inequality
\cite{Berezin1972,LiYau1983} gives
\begin{equation}\label{eq:li-yau-sum}
\sum_{j=1}^k\lambda_j
\geq \frac{n}{n+2}C_nV^{-2/n}k^{1+2/n}.
\end{equation}
Since $\sum_{j=1}^k\lambda_j\leq k\lambda_k$, it implies
\begin{equation}\label{eq:li-yau-individual}
\lambda_k\geq \frac{n}{n+2}C_nV^{-2/n}k^{2/n}.
\end{equation}
The leading constant in the sum inequality is asymptotically sharp, but
the deduction of \cref{eq:li-yau-individual} discards the spectral
deficit $k\lambda_k-\sum_{j=1}^k\lambda_j$.  This distinction motivates
the two results of the present paper.

The first result is a universal individual-eigenvalue estimate.  Exact
coordinate-commutator identities relate the spectral deficit to a
projected-gradient density.  At the same time, the scalar and vector
Fourier projections of a plane wave satisfy a nonlinear pointwise
compatibility condition.  After normalization, these relations yield
two densities $f$ and $g$ satisfying
\[
\sqrt t\,f\leq\sqrt{fg}+\sqrt{(1-f)(f-g)},
\qquad
\int_{\R^n}\left(tf+\frac2n g\right)\dd\eta
\leq\int_{\R^n}f\dd\eta.
\]
Eliminating $g$ pointwise reduces the problem to a one-dimensional
variational integral and gives
\begin{equation}\label{eq:intro-main-result}
\lambda_k\geq C_n^{\mathrm{new}}C_nV^{-2/n}k^{2/n},
\qquad
\frac{n}{n+2}<c_n^{\mathrm{env}}<C_n^{\mathrm{new}}<1.
\end{equation}
Here $n/(n+2)$ is the direct Li--Yau coefficient, while
$c_n^{\mathrm{env}}$ is the coefficient generated by the uncoupled
envelope
$F_k(\xi)\leq M_\Omega\min\{1,\lambda_k/|\xi|^2\}$,
$M_\Omega=V/(2\pi)^n$, studied in \cite{WangXie}.  We recall its
variational definition in \cref{eq:cenv} so that all constants are
self-contained.  The two strict inequalities in
\cref{eq:intro-main-result} separate the gain from the scalar energy cap
from the additional gain due to commutator coupling; see
\cref{prop:constant-comparison,tab:constants} for the proof and numerical
comparison.

The second result improves the eigenvalue sum at each finite index.  For
a bounded open set we introduce a positive trace-class
frequency-concentration operator $T_r$ on $L^2(B_r)$.  At the Li--Yau
filling radius $R_k$, one has $\Tr(T_{R_k})=k$, whereas analytic
uniqueness implies $\Tr(T_{R_k}^2)<k$.  The Ky Fan principle then forces
some of the Fourier mass of the first $k$ eigenfunctions to lie outside
$B_{R_k}$.  A two-region bathtub problem converts this leakage into
\begin{equation}\label{eq:intro-sum-result}
\sum_{j=1}^k\lambda_j
\geq \mathcal C_k^{(2)}(\Omega)
\frac{n}{n+2}C_nV^{-2/n}k^{1+2/n},
\qquad \mathcal C_k^{(2)}(\Omega)>1.
\end{equation}
The multiplier has an explicit squared-trace representation, decreases
strictly to one, and satisfies
\begin{equation}\label{eq:intro-quantitative-gap}
\mathcal C_k^{(2)}(\Omega)-1\geq\kappa_nk^{-2/n},
\qquad
\kappa_n=\frac{n+2}{2(n+1)^2}
\left(\frac{\omega_{n-1}}{\omega_n}\right)^2
\left(\frac34\right)^{n-1}.
\end{equation}
Equivalently, \cref{eq:intro-sum-result} contains the additional
volume-only term
\[
\frac{n}{n+2}\kappa_nC_nV^{-2/n}k.
\]
The decay of the relative correction is consistent with the sharp
semiclassical leading coefficient.

These two statements recover different pieces of information.  The
commutator argument strengthens a Weyl-scale lower bound for a single
eigenvalue, while the concentration argument strengthens the sharp sum
bound without changing its leading coefficient.  Neither result proves
P\'olya's conjecture: $C_n^{\mathrm{new}}<1$, and
$\mathcal C_k^{(2)}(\Omega)\to1$.  The role of the present estimates is
therefore to quantify two losses in the classical Fourier argument and
to identify where additional information would have to enter.

We briefly compare these conclusions with the closest available results.
For individual eigenvalues, the strict chain in
\cref{eq:intro-main-result} improves both the direct Li--Yau coefficient
and the scalar-envelope coefficient from our previous work \cite{WangXie}.
Its advantage is uniformity: it holds for every $k$ on every finite-measure open set.
By contrast, Jiang and Lin proved that for every
$\varepsilon\in(0,1)$ and every bounded Lipschitz domain,
\begin{equation}\label{eq:intro-jiang-lin}
\lambda_k\geq(1+\varepsilon)^{-2/n}
C_nV^{-2/n}k^{2/n}
\qquad\text{whenever }\lambda_k\geq
\Lambda_\varepsilon(\Omega),
\end{equation}
where the explicit threshold depends on $\varepsilon$ and the geometry
of $\Omega$ \cite{JiangLin2026}.  Thus
\cref{eq:intro-jiang-lin} approaches the P\'{o}lya coefficient at high
energy, whereas our fixed coefficient is valid without a threshold or
boundary regularity.  The full coefficient one is known for special
geometries, including Euclidean balls and sufficiently thin product
domains \cite{FLPSBalls2023,HeWang2026}; other universal individual
bounds of a different form appear in \cite{Steinerberger2024}.

For eigenvalue sums, write
\[
L_k^{\mathrm{LY}}
:=\frac{n}{n+2}C_nV^{-2/n}k^{1+2/n}.
\]
Melas proved the inertia-dependent correction
\begin{equation}\label{eq:intro-melas}
\sum_{j=1}^k\lambda_j
\geq L_k^{\mathrm{LY}}
+c_n^{\mathrm M}\frac{V}{I(\Omega)}k,
\qquad
I(\Omega)=\inf_{a\in\R^n}\int_\Omega|x-a|^2\dd x,
\end{equation}
with subsequent refinements in
\cite{Melas2003,Ilyin2010,GJL2025,JiLuo2026}.  Our quantitative estimate
instead gives
\begin{equation}\label{eq:intro-volume-linear-comparison}
\sum_{j=1}^k\lambda_j
\geq L_k^{\mathrm{LY}}
+\frac{n}{n+2}\kappa_nC_nV^{-2/n}k
\end{equation}
for every $k$ on every bounded open set.  Both remainders are linear in
$k$, but \cref{eq:intro-volume-linear-comparison} depends only on the
volume and not on $I(\Omega)$.

Several stronger-looking corrections belong to different regimes.
Volume-only strict improvements for order-one Riesz means under minimal
assumptions were obtained in \cite{FLP2025}, with an exponentially small
relative gain at high energy.  Boundary-sensitive estimates and two-term
Riesz-mean asymptotics can reach the larger second Weyl scale under
additional geometric, convexity, or Lipschitz hypotheses
\cite{KVW2009,Weidl2008,GLW2011,FrankLarson2025,FrankLarson2026}.
These results exploit information unavailable on an arbitrary bounded
open set.  The point of
\cref{eq:intro-sum-result,eq:intro-quantitative-gap} is instead an
explicit finite-$k$, volume-only improvement requiring no boundary
smoothness, convexity, connectedness, or moment-of-inertia term.

The paper is organized as follows.  \Cref{sec:preliminaries} records the
Fourier identities and the radial-envelope bathtub principle.
\Cref{sec:commutator-projection} establishes the commutator identity, the
coupled projection inequality, and the normalized slack formula.
\Cref{sec:main-results} defines $C_n^{\mathrm{new}}$, proves the
individual bound, and compares the three variational coefficients.
\Cref{sec:concentration} develops the concentration operator and its
trace formulas, and \cref{sec:strict-sum} proves the strict sum estimate
and its explicit excess.  The final section summarizes the scope and
limitations of the two methods in this paper.

\section{Preliminaries}\label{sec:preliminaries}

We begin by recording the spectral and Fourier identities on which the
argument rests. Let $\{\phi_j\}_{j\geq1}$ be a real-valued orthonormal basis of
Dirichlet eigenfunctions in $L^2(\Omega)$. We use complex Hilbert spaces
when pairing these functions with plane waves, with inner product linear
in its first argument. Thus
\begin{equation}\label{eq:eigenproblem}
-\Delta\phi_j=\lambda_j\phi_j,
\qquad \phi_j\in H_0^1(\Omega),
\qquad \int_\Omega\phi_i\phi_j\dd x=\delta_{ij}.
\end{equation}
We use the unitary Fourier transform
\begin{equation}\label{eq:fourier}
\wh u(\xi)=(2\pi)^{-\frac{n}{2}}
\int_{\R^n}u(x)e^{-\mathrm{i}x\cdot\xi}\dd x,
\end{equation}
where every function in $H_0^1(\Omega)$ is extended by zero to $\R^n$. The
Fourier density associated with the first $k$ eigenfunctions is
\begin{equation}\label{eq:density}
F_k(\xi)=\sum_{j=1}^k\abs{\wh\phi_j(\xi)}^2.
\end{equation}
Plancherel's theorem and the eigenvalue equation then give the two basic
identities
\begin{equation}\label{eq:mass-moment}
\int_{\R^n}F_k(\xi)\dd\xi=k,
\qquad
\int_{\R^n}\abs\xi^2F_k(\xi)\dd\xi
=\sum_{j=1}^k\lambda_j.
\end{equation}
Thus the spectral problem has been converted into a constrained moment problem:
the total mass of $F_k$ is known exactly, its second moment is the eigenvalue
sum, and the remaining task is to obtain the sharpest available pointwise
majorant.

Before deriving that majorant, we formulate the rearrangement principle that
will turn it into a lower bound for the second moment. The following statement
is a radial-obstacle version of the standard bathtub principle; see
\cite{WangXie} and Lieb and Loss~\cite{LiebLoss2001}.

\begin{lemma}[Bathtub principle under a radial envelope]\label{lem:bathtub}
Let $h:[0,\infty)\to[0,\infty)$ be measurable and nonincreasing, and suppose
that $h(\abs\cdot)\in L^1_{\mathrm{loc}}(\R^n)$. Let $0\leq M<\infty$ and
assume that
\[
\int_{\R^n}h(\abs\xi)\dd\xi\geq M.
\]
If a measurable function $f:\R^n\to[0,\infty)$ satisfies
\[
0\leq f(\xi)\leq h(\abs\xi)\quad\text{a.e.},
\qquad
\int_{\R^n}f(\xi)\dd\xi=M,
\]
then there exists $R\in[0,\infty]$ such that
\[
\int_{B_R}h(\abs\xi)\dd\xi=M.
\]
Consequently, if
\[
g_R(\xi):=h(\abs\xi)\ind_{B_R}(\xi),
\]
then
\begin{equation}\label{eq:bathtub}
\int_{\R^n}\abs\xi^2f(\xi)\dd\xi
\geq
\int_{B_R}\abs\xi^2h(\abs\xi)\dd\xi.
\end{equation}
The integrals in \cref{eq:bathtub} may take the value $+\infty$.
\end{lemma}

\begin{proof}
Set
\[
H(R)=\int_{B_R}h(|\xi|)\dd\xi.
\]
Local integrability implies that $H$ is continuous and nondecreasing,
with $H(0)=0$ and
$H(R)\uparrow\int_{\R^n}h(|\xi|)\dd\xi$ as $R\to\infty$.
Hence a radius $R\in[0,\infty]$ with $H(R)=M$ exists.  If
$0<R<\infty$, then $g_R=h(|\cdot|)\ind_{B_R}$ has the same mass as
$f$, and
\[
(|\xi|^2-R^2)\bigl(f(\xi)-g_R(\xi)\bigr)\geq0
\quad\text{for a.e. }\xi\in\R^n.
\]
Indeed, both factors are nonpositive on $B_R$, whereas both are
nonnegative on its complement.  Integration and cancellation of the
equal-mass terms multiplied by $R^2$ give \cref{eq:bathtub}.  If $R=0$,
then $M=0$ and the assertion is immediate.  If $R=\infty$, then
$M=\int h(|\xi|)\dd\xi$; the inequalities $0\leq f\leq h(|\cdot|)$
and equality of their integrals imply $f=h(|\cdot|)$ almost everywhere,
so equality holds in \cref{eq:bathtub}.
\end{proof}

The lemma isolates the classical variational part of the argument.  The
usual $L^2$ and energy estimates supply a scalar radial cap for $F_k$,
whereas the next section retains more information by coupling that density
to the low-mode gradient projection and measuring the integrated slack
through a commutator identity.

\section{Commutator Identities and Coupled Fourier Projections}
\label{sec:commutator-projection}

The preceding radial bathtub principle becomes substantially stronger when
the scalar Fourier density is not treated in isolation. We now retain its
interaction with the low-mode part of the gradient. For the commutator
calculation we first assume that $\Omega$ is bounded. The final eigenvalue
bound is extended to arbitrary finite-measure open sets by exhaustion.

\subsection{Coordinate commutators and an exact deficit identity}

For $1\leq\alpha\leq n$ and $i,j\geq1$, set
\begin{equation}\label{eq:a-d-matrices}
a_{ij}^{(\alpha)}=\int_\Omega x_\alpha\phi_i\phi_j\dd x,
\qquad
d_{ij}^{(\alpha)}=\int_\Omega\phi_i\,\partial_\alpha\phi_j\dd x.
\end{equation}
% The identities below are, respectively, the coordinate commutator identity
% and the Thomas--Reiche--Kuhn sum rule \cite{LevitinParnovski2002}.  They are
% the starting point for an exact finite-dimensional slack, which will later
% be converted into pointwise and integrated constraints on the Fourier
% density.
The first identity below is the eigenfunction matrix-element form of
the coordinate commutator;
see \cite[Eq.~(2.6) and Example~4.1]{LevitinParnovski2002}.
The second is the coordinate Thomas--Reiche--Kuhn sum rule
\cite[Eq.~(1.2)]{LevitinParnovski2002}.
We give a direct proof using the weak eigenvalue equation. 

\begin{lemma}[Coordinate commutator identities]\label{lem:commutator-identities}
For every $i,j\geq1$ and $1\leq\alpha\leq n$,
\begin{align}
(\lambda_i-\lambda_j)a_{ij}^{(\alpha)}
&=-2d_{ij}^{(\alpha)},\label{eq:gap-commutator}\\
\sum_{i=1}^\infty(\lambda_i-\lambda_j)
|a_{ij}^{(\alpha)}|^2&=1.\label{eq:trk}
\end{align}
The series in \cref{eq:trk} is absolutely convergent.
\end{lemma}

\begin{proof}
Because $\Omega$ is bounded, $x_\alpha\in W^{1,\infty}(\Omega)$, and
multiplication by either $x_\alpha$ or $x_\alpha^2$ preserves
$H_0^1(\Omega)$.  The weak eigenvalue equation for $\phi_i$, tested with
$x_\alpha\phi_j$, and the equation for $\phi_j$, tested with
$x_\alpha\phi_i$, give
\begin{align}
\lambda_i a_{ij}^{(\alpha)}
&=\int_\Omega x_\alpha\nabla\phi_i\cdot\nabla\phi_j\dd x
  +\int_\Omega\phi_j\partial_\alpha\phi_i\dd x,\label{equality_1}\\
\lambda_j a_{ij}^{(\alpha)}
&=\int_\Omega x_\alpha\nabla\phi_i\cdot\nabla\phi_j\dd x
  +\int_\Omega\phi_i\partial_\alpha\phi_j\dd x.\label{equality_2}
\end{align}
Here $\phi_i\phi_j\in W_0^{1,1}(\Omega)$, so weak integration by parts
implies
\begin{align}
\int_\Omega\phi_j\partial_\alpha\phi_i\dd x
=-\int_\Omega\phi_i\partial_\alpha\phi_j\dd x
=-d_{ij}^{(\alpha)}.\label{equality_3}
\end{align}
With \eqref{equality_3} and subtracting the preceding two identities \eqref{equality_1} and \eqref{equality_2} therefore proves 
\cref{eq:gap-commutator}.  Notice also that
$a_{ji}^{(\alpha)}=a_{ij}^{(\alpha)}$ and
$d_{ji}^{(\alpha)}=-d_{ij}^{(\alpha)}$; in particular,
$d_{ii}^{(\alpha)}=0$.

For the sum rule, put $u=x_\alpha\phi_j\in H_0^1(\Omega)$.  The spectral
representation of the Dirichlet form $x_\alpha\phi_j=\sum_{i=1}^\infty a_{ij}^{(\alpha)}\phi_i$ 
and then  Parseval's identity yield 
\begin{align}
\sum_{i=1}^{\infty}\lambda_i|a_{ij}^{(\alpha)}|^2
&=\int_\Omega|\nabla u|^2\dd x,\label{equality_4}\\
\sum_{i=1}^{\infty}|a_{ij}^{(\alpha)}|^2
&=\int_\Omega x_\alpha^2\phi_j^2\dd x.\label{equality_5}
\end{align}
Expanding the first integral makes the constant in the sum rule explicit:
\begin{align}
\int_\Omega|\nabla(x_\alpha\phi_j)|^2\dd x
=\int_\Omega x_\alpha^2|\nabla\phi_j|^2\dd x
 +2\int_\Omega x_\alpha\phi_j\partial_\alpha\phi_j\dd x
 +\int_\Omega\phi_j^2\dd x.\label{equality_6}
\end{align}
On the other hand, testing the $j$-th eigenvalue equation with
$x_\alpha^2\phi_j$ gives
\begin{align}
\lambda_j\int_\Omega x_\alpha^2\phi_j^2\dd x
=\int_\Omega x_\alpha^2|\nabla\phi_j|^2\dd x
 +2\int_\Omega x_\alpha\phi_j\partial_\alpha\phi_j\dd x.\label{equality_7}
\end{align}
After subtraction, 
The combination of \eqref{equality_4}, \eqref{equality_5}, \eqref{equality_6}, \eqref{equality_7} and the $L^2$ normalization
$\int_\Omega\phi_j^2\dd x=1$ proves \cref{eq:trk}.  Finally,
\[
\sum_{i=1}^{\infty}|\lambda_i-\lambda_j|
 |a_{ij}^{(\alpha)}|^2
\leq
\sum_{i=1}^{\infty}(\lambda_i+\lambda_j)
 |a_{ij}^{(\alpha)}|^2<\infty,
\]
by the two spectral identities above.  Thus the series in
\cref{eq:trk} is absolutely convergent.
\end{proof}

Fix $k\geq1$, write
\begin{equation}\label{eq:new-deficit-notation}
\Lambda=\lambda_k,
\qquad S_k=\sum_{j=1}^k\lambda_j,
\qquad D_k=k\Lambda-S_k,
\end{equation}
and introduce
\begin{align}
I_k&=\sum_{\alpha=1}^n\sum_{i=1}^k\sum_{j=1}^k
|d_{ij}^{(\alpha)}|^2,\label{eq:Ik-definition}\\
\mathcal B_k&=\sum_{\alpha=1}^n\sum_{j=1}^k\sum_{i>k}
(\Lambda-\lambda_j)(\lambda_i-\lambda_j)
|a_{ij}^{(\alpha)}|^2.\label{eq:Bcal-definition}
\end{align}
Every term in $\mathcal B_k$ is nonnegative, so $\mathcal B_k\geq0$.

The following identity is the Dirichlet Laplacian specialization of
the trace formula of Harrell and Stubbe
\cite[Theorem~1, Eq.~(4)]{HarrellStubbe1997}.
We give a direct proof in our notation, which also covers the case
where the cutoff at $k$ splits a multiple eigenvalue.
\begin{proposition}[Exact commutator slack]\label{prop:exact-commutator-slack}
The quantities in \cref{eq:new-deficit-notation,eq:Ik-definition,eq:Bcal-definition}
satisfy
\begin{equation}\label{eq:exact-commutator-slack}
nD_k=2I_k+\mathcal B_k.
\end{equation}
In particular,
\begin{equation}\label{eq:Ik-Dk}
I_k\leq\frac n2D_k.
\end{equation}
\end{proposition}

\begin{proof}
Multiply \cref{eq:trk} by $\Lambda-\lambda_j$ and sum over
$1\leq j\leq k$ and $1\leq\alpha\leq n$.  The absolute convergence
established in \cref{lem:commutator-identities} justifies splitting each
$i$-sum at $k$, and hence
\begin{align*}
nD_k
&=\sum_{\alpha=1}^n\sum_{j=1}^k(\Lambda-\lambda_j)
  \sum_{i=1}^{\infty}(\lambda_i-\lambda_j)
  |a_{ij}^{(\alpha)}|^2\\
&=\sum_{\alpha=1}^n\sum_{i=1}^k\sum_{j=1}^k
  (\Lambda-\lambda_j)(\lambda_i-\lambda_j)
  |a_{ij}^{(\alpha)}|^2+\mathcal B_k.
\end{align*}
It remains to identify the finite double sum.  Since
$a_{ji}^{(\alpha)}=a_{ij}^{(\alpha)}$, symmetrization in $i$ and $j$
gives, for each $\alpha$,
\begin{align*}
&\sum_{i=1}^k\sum_{j=1}^k(\Lambda-\lambda_j)(\lambda_i-\lambda_j)
  |a_{ij}^{(\alpha)}|^2\\
&\quad=\frac12\sum_{i=1}^k\sum_{j=1}^k
 \Bigl[(\Lambda-\lambda_j)(\lambda_i-\lambda_j)
 +(\Lambda-\lambda_i)(\lambda_j-\lambda_i)\Bigr]
 |a_{ij}^{(\alpha)}|^2\\
&\quad=\frac12\sum_{i=1}^k\sum_{j=1}^k
 (\lambda_i-\lambda_j)^2|a_{ij}^{(\alpha)}|^2
 =2\sum_{i=1}^k\sum_{j=1}^k|d_{ij}^{(\alpha)}|^2,
\end{align*}
where the last equality follows from \cref{eq:gap-commutator}.  Summing
over $\alpha$ proves \cref{eq:exact-commutator-slack}.  Since
$\mathcal B_k\geq0$, \cref{eq:Ik-Dk} follows immediately.  In
\cref{prop:normalized-moments}, the same identity will reappear as the
exact slack in the normalized integral constraint.
\end{proof}

\subsection{The coupled scalar--vector projection inequality}\label{Section_3.2}

Set
\begin{equation}\label{eq:M-Omega-definition}
M_\Omega=\frac{V}{(2\pi)^n}
\end{equation}
and let $P_k$ be the orthogonal projection in $L^2(\Omega)$ onto
$E_k=\operatorname{span}\{\phi_1,\ldots,\phi_k\}$. Its componentwise
extension to $L^2(\Omega;\C^n)$ is denoted by $\mathbf P_k$. For
$\xi\in\R^n$, define
\[
E_\xi(x)=(2\pi)^{-n/2}e^{-\ii x\cdot\xi},
\qquad q_\xi=P_kE_\xi,
\]
and retain the notation $F_k(\xi)=\|q_\xi\|_2^2$. We also set
\begin{equation}\label{eq:T-G-definitions}
T_k(\xi)=\|\nabla q_\xi\|_2^2,
\qquad
G_k(\xi)=\|\mathbf P_k\nabla q_\xi\|_2^2.
\end{equation}
If $c(\xi)=(c_j(\xi))_{1\leq j\leq k}=(\widehat\phi_j(\xi))_{1\leq j\leq k}$ and
$D^{(\alpha)}=(d_{ij}^{(\alpha)})_{1\leq i,j\leq k}$, then
\begin{equation}\label{eq:G-matrix-form}
G_k(\xi)=\sum_{\alpha=1}^n
\left\|D^{(\alpha)}c(\xi)\right\|_{\C^k}^2.
\end{equation}
Consequently, Plancherel's identity gives the exact integral relation
\begin{equation}\label{eq:G-integral}
\int_{\R^n}G_k(\xi)\dd\xi=I_k.
\end{equation}
Moreover,
\begin{equation}\label{eq:F-T-G-basic}
0\leq G_k(\xi)\leq T_k(\xi)
=\sum_{j=1}^k\lambda_j|\widehat\phi_j(\xi)|^2
\leq\Lambda F_k(\xi).
\end{equation}

\begin{proposition}[Coupled Fourier projection]\label{prop:coupled-projection}
For every $\xi\neq0$,
\begin{equation}\label{eq:coupled-projection}
|\xi|F_k(\xi)
\leq\sqrt{F_k(\xi)G_k(\xi)}
+\sqrt{\bigl(M_\Omega-F_k(\xi)\bigr)
  \bigl(T_k(\xi)-G_k(\xi)\bigr)}.
\end{equation}
In particular,
\begin{equation}\label{eq:coupled-projection-Lambda}
|\xi|F_k(\xi)
\leq\sqrt{F_k(\xi)G_k(\xi)}
+\sqrt{\bigl(M_\Omega-F_k(\xi)\bigr)
  \bigl(\Lambda F_k(\xi)-G_k(\xi)\bigr)}.
\end{equation}
\end{proposition}

\begin{proof}
Fix $\xi\neq0$ and set
\[
Y_\xi=-\ii\frac{\xi}{|\xi|}E_\xi
\in L^2(\Omega;\C^n).
\]
Since $|E_\xi|=(2\pi)^{-n/2}$, the definition
\cref{eq:M-Omega-definition} gives
\begin{equation}\label{eq:Y-norm}
\|Y_\xi\|_2^2=M_\Omega.
\end{equation}
The componentwise projection acts only on the scalar factor of
$Y_\xi$, and therefore
\begin{equation}\label{eq:PY-formula}
\mathbf P_kY_\xi
=-\ii\frac{\xi}{|\xi|}P_kE_\xi
=-\ii\frac{\xi}{|\xi|}q_\xi,
\qquad
\|\mathbf P_kY_\xi\|_2^2=\|q_\xi\|_2^2=F_k(\xi).
\end{equation}

We next compute the relevant pairing. Because
$q_\xi\in E_k\subset H_0^1(\Omega)$, weak integration by parts is
legitimate. Using
$\operatorname{div}\overline{Y_\xi}=-|\xi|\overline{E_\xi}$, we obtain
\begin{align}
\langle\nabla q_\xi,Y_\xi\rangle
&=-\int_\Omega q_\xi\,
\operatorname{div}\overline{Y_\xi}\dd x\notag\\
&=|\xi|\langle q_\xi,E_\xi\rangle
=|\xi|\langle P_kE_\xi,E_\xi\rangle\notag\\
&=|\xi|\|P_kE_\xi\|_2^2
=|\xi|F_k(\xi).
\label{eq:projection-pairing}
\end{align}
Here the penultimate equality uses that $P_k$ is an orthogonal
projection.

Decompose both vector fields with respect to the orthogonal splitting
\[
L^2(\Omega;\C^n)
=\Ran\mathbf P_k\oplus\Ker\mathbf P_k.
\]
The mixed pairings vanish, so
\begin{align}
\langle\nabla q_\xi,Y_\xi\rangle
=\langle\mathbf P_k\nabla q_\xi,
  \mathbf P_kY_\xi\rangle+
\langle(I-\mathbf P_k)\nabla q_\xi,
  (I-\mathbf P_k)Y_\xi\rangle.
\label{eq:orthogonal-splitting-pairing}
\end{align}
By \cref{eq:T-G-definitions,eq:Y-norm,eq:PY-formula} and orthogonality,
\begin{align*}
\|\mathbf P_k\nabla q_\xi\|_2^2&=G_k(\xi),
&\|(I-\mathbf P_k)\nabla q_\xi\|_2^2
&=T_k(\xi)-G_k(\xi),\\
\|\mathbf P_kY_\xi\|_2^2&=F_k(\xi),
&\|(I-\mathbf P_k)Y_\xi\|_2^2
&=M_\Omega-F_k(\xi).
\end{align*}
Applying these identities, the triangle inequality and Cauchy--Schwarz to the two
orthogonal components in \cref{eq:orthogonal-splitting-pairing}, and then using
\cref{eq:projection-pairing}, proves \cref{eq:coupled-projection}.
Finally, \cref{eq:F-T-G-basic} gives
$T_k(\xi)-G_k(\xi)\leq\Lambda F_k(\xi)-G_k(\xi)$.
Substitution into \cref{eq:coupled-projection} proves
\cref{eq:coupled-projection-Lambda}.
\end{proof}

The scalar consequences include both familiar Bessel bounds.
\begin{corollary}[The uncoupled envelope]\label{prop:envelope}
For almost every $\xi\in\R^n$,
\begin{equation}\label{eq:l2-bessel}
F_k(\xi)\leq M_\Omega,
\end{equation}
and, for $\xi\neq0$,
\begin{equation}\label{eq:projection-envelope}
F_k(\xi)\leq M_\Omega\min\left\{1,\frac{\Lambda}{|\xi|^2}\right\}.
\end{equation}
\end{corollary}

\begin{proof}
The first estimate is the scalar Bessel inequality. Applying
Cauchy--Schwarz directly to \cref{eq:coupled-projection}, or squaring the
projection identity before splitting, gives
$|\xi|^2F_k(\xi)^2\leq M_\Omega T_k(\xi)
\leq M_\Omega\Lambda F_k(\xi)$.
Combining these two estimates proves the assertion.
\end{proof}
The scalar envelope \cref{eq:projection-envelope} was also used in
\cite{WangXie}.  Its associated variational coefficient is denoted here
by $c_n^{\mathrm{env}}$ in \cref{eq:intro-main-result}
and is defined independently in \cref{eq:cenv}.

\subsection{Normalized constraints and the exact integrated slack}

Put $p=n/2$, rescale $\xi=\sqrt\Lambda\,\eta$, and write
\begin{equation}\label{eq:f-g-normalization}
t=|\eta|^2,
\qquad
f(\eta)=\frac{F_k(\sqrt\Lambda\,\eta)}{M_\Omega},
\qquad
g(\eta)=\frac{G_k(\sqrt\Lambda\,\eta)}{\Lambda M_\Omega}.
\end{equation}
For $t\geq0$ define
\begin{equation}\label{eq:m-Psi-definitions}
m(t)=\begin{cases}
1,&t\in [0,1],\\
t^{-1},&t>1,
\end{cases}
\qquad
\Psi(t,u)=u\left[
\sqrt t\,u-\sqrt{(1-u)(1-tu)}
\right]_+^2.
\end{equation}

\begin{lemma}[Basic properties of the penalty]\label{lem:Psi-properties}
On its natural domain $0\leq u\leq m(t)$, the function $\Psi$ is
continuous, nonnegative, and nondecreasing separately in $t$ and $u$.
For the monotonicity in $t$, this means that $t_1\leq t_2$ and
$u\leq m(t_2)$, so that both $(t_1,u)$ and $(t_2,u)$ lie in the natural
domain.
Moreover,
\begin{equation}\label{eq:Psi-zero-set}
\Psi(t,u)=0\quad\Longleftrightarrow\quad
u\leq\frac1{1+t},
\end{equation}
and, for $t>0$,
\begin{equation}\label{eq:Psi-at-cap}
\Psi(t,m(t))=
\begin{cases}
t,&0<t\leq1,\\
t^{-2},&t\geq1.
\end{cases}
\end{equation}
\end{lemma}

\begin{proof}
Continuity and nonnegativity are immediate.  Put
\[
B(t,u)=\sqrt t\,u-\sqrt{(1-u)(1-tu)}.
\]
For fixed $u$, the first term is nondecreasing in $t$ and the expression
under the second square root is nonincreasing.  Hence $B$ and
$u[B]_+^2$ are nondecreasing in $t$.  For fixed $t$, note that
\[
\frac{\partial}{\partial u}\bigl[(1-u)(1-tu)\bigr]
=-(1+t)+2tu\leq0
\]
on $0\leq u\leq m(t)$: if $t\leq1$, then $2tu\leq2t\leq1+t$, while if
$t\geq1$, then $2tu\leq2\leq1+t$.  Thus $B$ is nondecreasing in $u$;
since both $u$ and $[B]_+^2$ are nonnegative and nondecreasing, so is
$\Psi(t,u)=u[B(t,u)]_+^2$.  Finally, the bracket is nonpositive exactly
when
\[
tu^2\leq(1-u)(1-tu),
\]
which reduces to $(1+t)u\leq1$ and proves
\cref{eq:Psi-zero-set}. Substitution of $m(t)=1$ for $t\leq1$ and
$m(t)=t^{-1}$ for $t\geq1$ proves \cref{eq:Psi-at-cap}.
\end{proof}

\begin{lemma}[Pointwise elimination of the gradient variable]
\label{lem:pointwise-elimination}
The normalized densities satisfy
\begin{equation}\label{eq:normalized-pointwise-system}
0\leq f\leq m(t),
\qquad 0\leq g\leq f,
\qquad
\sqrt t\,f\leq\sqrt{fg}+\sqrt{(1-f)(f-g)},
\end{equation}
and hence
\begin{equation}\label{eq:g-Psi}
g(\eta)\geq\Psi(t,f(\eta)).
\end{equation}
\end{lemma}

\begin{proof}
By \cref{eq:projection-envelope}, evaluated at
$\xi=\sqrt\Lambda\,\eta$,
\[
0\leq f(\eta)\leq
\min\left\{1,\frac1{|\eta|^2}\right\}=m(t),
\]
where the value at $t=0$ is understood through the first branch of
$m$.  Likewise, \cref{eq:F-T-G-basic} gives $0\leq g\leq f$.  Finally,
substituting $F_k=M_\Omega f$, $G_k=\Lambda M_\Omega g$, and
$|\xi|=\sqrt\Lambda\sqrt t$ into
\cref{eq:coupled-projection-Lambda}, and then dividing by
$\sqrt\Lambda M_\Omega$, yields
\[
\sqrt t\,f\leq\sqrt{fg}+\sqrt{(1-f)(f-g)}.
\]
At $t=0$ this inequality is automatic.  Thus
\cref{eq:normalized-pointwise-system} holds.

It remains to eliminate $g$.  The bound $f\leq m(t)$ implies both
$0\leq f\leq1$ and $0\leq tf\leq1$, so every square root below is real.
If $f>0$, set
\[
s=\sqrt{\frac gf},\qquad
a=\sqrt f,\qquad b=\sqrt{1-f},\qquad r=\sqrt{tf}.
\]
Then $a^2+b^2=1$, $0\leq r,s\leq1$, and division of the preceding
inequality by $\sqrt f$ gives the scalar constraint
\[
r\leq as+b\sqrt{1-s^2}.
\]
Write $a=\sin\alpha$ and $s=\sin\theta$ with
$\alpha,\theta\in[0,\pi/2]$.  Its right-hand side is
\[
\sin\alpha\sin\theta+\cos\alpha\cos\theta
=\cos(\theta-\alpha).
\]
If $r\leq b$, then $s=0$ is admissible.  If $r>b$, the first admissible
angle is $\theta=\alpha-\arccos r$, and hence
\[
s\geq ar-b\sqrt{1-r^2}.
\]
Both cases can therefore be written as
\[
s\geq
\left[ar-b\sqrt{1-r^2}\right]_+
=\left[
\sqrt t\,f-\sqrt{(1-f)(1-tf)}
\right]_+.
\]
Since $g=fs^2$, squaring and multiplying by $f$ proves
$g\geq\Psi(t,f)$.  If $f=0$, the already proved inequality
$0\leq g\leq f$ gives $g=0=\Psi(t,0)$.  The inactive branch in this
argument is exactly the zero set described in \cref{eq:Psi-zero-set}.
\end{proof}

\begin{proposition}[Moments and exact slack]\label{prop:normalized-moments}
The zeroth and second moments and the projected-gradient mass are
\begin{align}
\int_{\R^n}f\dd\eta
&=\frac{k}{M_\Omega\Lambda^p},\label{eq:f-mass}\\
\int_{\R^n}t f\dd\eta
&=\frac{S_k}{M_\Omega\Lambda^{p+1}},\label{eq:f-moment}\\
\int_{\R^n}g\dd\eta
&=\frac{I_k}{M_\Omega\Lambda^{p+1}}.\label{eq:g-mass}
\end{align}
Moreover,
\begin{equation}\label{eq:coupled-integrated-constraint}
\int_{\R^n}\left(tf+\frac1p g\right)\dd\eta
\leq\int_{\R^n}f\dd\eta,
\end{equation}
and the slack is exactly
\begin{equation}\label{eq:normalized-exact-slack}
\int_{\R^n}f\dd\eta
-\int_{\R^n}\left(tf+\frac1p g\right)\dd\eta
=\frac{\mathcal B_k}{nM_\Omega\Lambda^{p+1}}.
\end{equation}
\end{proposition}

\begin{proof}
Under the change of variables $\xi=\sqrt\Lambda\,\eta$, one has
$\dd\eta=\Lambda^{-n/2}\dd\xi=\Lambda^{-p}\dd\xi$ and
$t=|\eta|^2=|\xi|^2/\Lambda$.  Consequently,
\begin{align*}
\int_{\R^n}f(\eta)\dd\eta
&=\frac{\Lambda^{-p}}{M_\Omega}
  \int_{\R^n}F_k(\xi)\dd\xi
 =\frac{k}{M_\Omega\Lambda^p},\\
\int_{\R^n}t f(\eta)\dd\eta
&=\frac{\Lambda^{-p-1}}{M_\Omega}
  \int_{\R^n}|\xi|^2F_k(\xi)\dd\xi
 =\frac{S_k}{M_\Omega\Lambda^{p+1}},\\
\int_{\R^n}g(\eta)\dd\eta
&=\frac{\Lambda^{-p-1}}{M_\Omega}
  \int_{\R^n}G_k(\xi)\dd\xi
 =\frac{I_k}{M_\Omega\Lambda^{p+1}}.
\end{align*}
Here the first two final equalities use \cref{eq:mass-moment}, and the
third uses \cref{eq:G-integral}.  This proves
\cref{eq:f-mass,eq:f-moment,eq:g-mass} and also exposes all scaling
factors in the normalization.

Using $p=n/2$ and the definition of $D_k$, the difference between the two
sides of \cref{eq:coupled-integrated-constraint} is
\begin{align*}
&\int_{\R^n}f\dd\eta
-\int_{\R^n}\left(tf+\frac1p g\right)\dd\eta\\
&\quad=\frac{k\Lambda-S_k-p^{-1}I_k}
 {M_\Omega\Lambda^{p+1}}
=\frac{D_k-\frac2n I_k}{M_\Omega\Lambda^{p+1}}
=\frac{\mathcal B_k}{nM_\Omega\Lambda^{p+1}},
\end{align*}
where the last equality is precisely
\cref{eq:exact-commutator-slack}.  This proves
\cref{eq:normalized-exact-slack}; its nonnegativity, already recorded
after \cref{eq:Bcal-definition}, gives
\cref{eq:coupled-integrated-constraint}.
\end{proof}

\begin{remark}[What is gained over the scalar envelope]
The cap \cref{eq:projection-envelope} remembers only $f\leq m(t)$.
The new system additionally forces $g\geq\Psi(t,f)$ and charges the
integral of $g$ through \cref{eq:coupled-integrated-constraint}. It is
this simultaneous pointwise and integrated restriction, rather than a
new scalar cap alone, that improves the principal coefficient below.
\end{remark}

\section{A Universal Commutator-Coupled Eigenvalue Constant}
\label{sec:main-results}

We now solve the dimensionless variational problem created by the coupled
constraints. For $a>1$, define
\begin{align}
L_{n,a}(t)
&=\max_{0\leq u\leq m(t)}
\left\{(a-t)u-\frac1p\Psi(t,u)\right\},\label{eq:Lna}\\
J_n(a)&=p\int_0^aL_{n,a}(t)t^{p-1}\dd t,\label{eq:Jn}\\
C_n^{\mathrm{new}}
&=\sup_{a>1}\left(\frac{a-1}{J_n(a)}\right)^{1/p}.
\label{eq:Cnew}
\end{align}
The maximization in \cref{eq:Lna} is one-dimensional and takes place on
a compact interval.  Moreover, $L_{n,a}(t)=0$ for $t\geq a$, because
the objective is then nonpositive and vanishes at $u=0$; for $0<t<a$
it is positive, since one may choose
$0<u\leq\min\{m(t),(1+t)^{-1}\}$, for which $\Psi(t,u)=0$.
Consequently $0<J_n(a)<\infty$, and $J_n(a)$ is directly computable
without any domain-dependent quantity.

\begin{theorem}[Commutator-coupled individual bound]
\label{thm:commutator-coupled-main}
Let $\Omega\subset\R^n$, $n\geq2$, be an open set with
$0<V=|\Omega|<\infty$. Then, for every $k\geq1$,
\begin{equation}\label{eq:commutator-coupled-main}
\lambda_k(\Omega)
\geq C_n^{\mathrm{new}}C_nV^{-2/n}k^{2/n}.
\end{equation}
\end{theorem}

\begin{proof}
Assume first that $\Omega$ is bounded and use the functions in
\cref{eq:f-g-normalization}. For every $a>1$,
\cref{eq:coupled-integrated-constraint} gives
\begin{align*}
(a-1)\int_{\R^n}f\dd\eta
&\leq\int_{\R^n}\left((a-t)f-\frac1p g\right)\dd\eta\\
&\leq\int_{\R^n}L_{n,a}(|\eta|^2)\dd\eta
=\omega_nJ_n(a),
\end{align*}
where \cref{lem:pointwise-elimination} was used in the second step.
By \cref{eq:f-mass},
\[
\frac{k}{M_\Omega\Lambda^p}\leq\frac{\omega_nJ_n(a)}{a-1}.
\]
If
\[
W_k=4\pi^2\left(\frac{k}{\omega_nV}\right)^{2/n},
\]
then $W_k^p=k/(M_\Omega\omega_n)$, and therefore
\[
\frac{\Lambda}{W_k}
\geq\left(\frac{a-1}{J_n(a)}\right)^{1/p}.
\]
Taking the supremum over $a>1$ proves the result for bounded domains.

For a general finite-measure open set, let
$\Omega_R=\Omega\cap B_R$. The spaces $H_0^1(\Omega_R)$ increase to a
form core of $H_0^1(\Omega)$, so the min--max principle gives
$\lambda_k(\Omega_R)\downarrow\lambda_k(\Omega)$, while
$|\Omega_R|\uparrow V$. Applying the bounded-domain result and passing
to the limit completes the proof.
\end{proof}

\subsection{Strict comparison with the uncoupled two-envelope method}

If the nonlinear penalty in \cref{eq:Lna} is discarded, one obtains
\begin{align}
J_{n,0}(a)&=p\int_0^a(a-t)m(t)t^{p-1}\dd t,\label{eq:Jn0}\\
c_n^{\mathrm{env}}
&=\max_{a>1}\left(\frac{a-1}{J_{n,0}(a)}\right)^{1/p}.
\label{eq:cenv}
\end{align}
A direct integration gives \cite{WangXie}
\begin{equation}\label{eq:Jn0-explicit}
J_{n,0}(a)=
\begin{cases}
a\log a+\dfrac12,&p=1,\\[5pt]
\dfrac{a^p-a}{p-1}+\dfrac1{p+1},&p>1.
\end{cases}
\end{equation}
This is the coefficient generated by the uncoupled scalar envelope
\cref{eq:projection-envelope}, which was analyzed in \cite{WangXie}.

For comparison, define the constant-density relaxation
\begin{equation}\label{eq:JnL2}
J_{n,L^2}(a)
=p\int_0^a(a-t)t^{p-1}\dd t
=\frac{a^{p+1}}{p+1}.
\end{equation}
It is the variational integral corresponding to the constant cap used in
the Li--Yau argument \cite{LiYau1983}.

\begin{proposition}[Strict improvement and P\'{o}lya obstruction]
\label{prop:constant-comparison}
For every $n\geq2$,
\begin{equation}\label{eq:strict-new-comparison}
\frac{n}{n+2}<c_n^{\mathrm{env}}<C_n^{\mathrm{new}}<1.
\end{equation}
\end{proposition}

\begin{proof}
Put
\[
\Phi_a(t,u)=(a-t)u-\frac1p\Psi(t,u),
\qquad 0\leq u\leq m(t).
\]
We prove the three comparisons separately.

\smallskip
\noindent\emph{Step 1: the coupled coefficient is strictly larger than
the scalar-envelope coefficient.}
Fix $a>1$ and $0<t<a$.  Since $a-t>0$ and $\Psi\geq0$,
\begin{equation}\label{eq:Lna-relaxed-comparison}
\Phi_a(t,u)\leq(a-t)u\leq(a-t)m(t).
\end{equation}
The admissible interval is compact and $\Phi_a(t,\cdot)$ is continuous,
so the maximum defining $L_{n,a}(t)$ is attained.  Equality in both
inequalities in \cref{eq:Lna-relaxed-comparison} would force
$u=m(t)$ and $\Psi(t,m(t))=0$.  This is impossible because
\cref{eq:Psi-at-cap} gives $\Psi(t,m(t))>0$ for $t>0$.  Therefore
\[
L_{n,a}(t)<(a-t)m(t),\qquad 0<t<a.
\]
Multiplication by $pt^{p-1}$ and integration yields
\begin{equation}\label{eq:J-strict-below-J0}
J_n(a)<J_{n,0}(a),\qquad a>1.
\end{equation}
The explicit formula \cref{eq:Jn0-explicit} shows that
\[
r_0(a):=\left(\frac{a-1}{J_{n,0}(a)}\right)^{1/p}
\]
is continuous on $(1,\infty)$ and tends to zero as $a\downarrow1$ and
as $a\to\infty$.  Hence it attains its maximum at some
$a_{\mathrm{env}}>1$.  By \cref{eq:J-strict-below-J0},
\[
C_n^{\mathrm{new}}
\geq\left(\frac{a_{\mathrm{env}}-1}
{J_n(a_{\mathrm{env}})}\right)^{1/p}
>r_0(a_{\mathrm{env}})=c_n^{\mathrm{env}}.
\]

\smallskip
\noindent\emph{Step 2: the scalar envelope strictly improves the direct
Li--Yau coefficient.}
Because $m(t)\leq1$ and $m(t)<1$ for $t>1$, for every $a>1$ we have
\begin{align*}
J_{n,L^2}(a)-J_{n,0}(a)
&=p\int_1^a(a-t)\bigl(1-t^{-1}\bigr)t^{p-1}\dd t>0.
\end{align*}
The function associated with the constant-density relaxation is
maximized at $a_L=1+1/p$, and
\begin{equation}\label{eq:L2-relaxation-maximum}
\max_{a>1}
\left(\frac{a-1}{J_{n,L^2}(a)}\right)^{1/p}
=\left(\frac{(p+1)(a_L-1)}{a_L^{p+1}}\right)^{1/p}
=\frac p{p+1}=\frac n{n+2}.
\end{equation}
Evaluating the scalar-envelope quotient at $a_L$ and using the strict
integral inequality gives
$c_n^{\mathrm{env}}>n/(n+2)$.

\smallskip
\noindent\emph{Step 3: the coupled relaxation still has coefficient
strictly below one.}
Define
\begin{equation}\label{eq:Ap-Tp}
A_p(T)=p\int_0^T\frac{t^{p-1}}{1+t}\dd t,
\qquad
T_p^p=2A_p(T_p),
\end{equation}
where the second equation has a unique solution $T_p>1$.  To see this,
let
$B_p(T)=T^p-2A_p(T)$,
so that
\begin{equation}\label{eq:Bp-derivative}
B_p'(T)=pT^{p-1}\frac{T-1}{T+1},
\end{equation}
for $T>0$.
Moreover, $A_p(1)>1/2$, hence $B_p(1)<0$.  If $p=1$, then
$A_1(T)=\log(1+T)$; if $p>1$, then
$A_p(T)=O(T^{p-1})$.  In either case $B_p(T)\to\infty$ as
$T\to\infty$.  Thus \cref{eq:Bp-derivative}
gives exactly one zero in $(1,\infty)$.

Consider the radial pair
\[
f_*(\eta)=\frac{\one_{\{|\eta|^2\leq T_p\}}}{1+|\eta|^2},
\qquad g_*(\eta)=0,
\]
and recall the radial integration formula
\begin{equation}\label{eq:radial-t-integration}
\int_{\R^n}h(|\eta|^2)\dd\eta
=p\omega_n\int_0^\infty h(t)t^{p-1}\dd t.
\end{equation}
For $0\leq t\leq T_p$, the identity
$f_*=1/(1+t)$ gives
\[
\sqrt t\,f_*=\sqrt{(1-f_*)f_*};
\]
outside this ball both sides vanish.  Also $0\leq f_*\leq m(t)$.
Consequently $(f_*,g_*)$ satisfies the complete pointwise system in
\cref{eq:normalized-pointwise-system}.  From
\cref{eq:radial-t-integration} and $T_p^p=2A_p(T_p)$,
\begin{align*}
\int_{\R^n}f_*\dd\eta
&=\omega_nA_p(T_p),\\
\int_{\R^n}|\eta|^2f_*\dd\eta
&=\omega_n\bigl(T_p^p-A_p(T_p)\bigr)
=\omega_nA_p(T_p).
\end{align*}
Thus the integrated constraint also holds, with equality.

Since $L_{n,a}(t)$ is the pointwise maximum of $\Phi_a(t,u)$ and is
zero for $t\geq a$, \cref{eq:radial-t-integration} gives, for every
$a>1$,
\[
\omega_nJ_n(a)\geq
\int_{\R^n}(a-|\eta|^2)f_*(\eta)\dd\eta
=(a-1)\omega_nA_p(T_p).
\]
It follows that
\begin{equation}\label{eq:Cflat}
C_n^{\mathrm{new}}
\leq C_n^{\mathrm{flat}}
:=A_p(T_p)^{-1/p}=\frac{2^{1/p}}{T_p}.
\end{equation}
It remains to verify that the last number is strictly smaller than one.
Let $T_0=2^{1/p}$.  After setting $t=T_0r$, Jensen's inequality for the
strictly convex function $r\mapsto(1+T_0r)^{-1}$ and the probability
density $pr^{p-1}$ on $(0,1)$ gives
\[
A_p(T_0)=2\int_0^1\frac{pr^{p-1}}{1+T_0r}\dd r
>\frac{2}{1+T_0p/(p+1)}\geq1,
\]
where the last inequality is equivalent to
$2^{1/p}\leq1+1/p$, which follows from Bernoulli's inequality
$(1+1/p)^p\geq2$.  Hence $B_p(T_0)=2-2A_p(T_0)<0$.
Since $B_p$ is strictly increasing on $(1,\infty)$ and vanishes at
$T_p$, we have $T_p>T_0$, and therefore
$C_n^{\mathrm{new}}\leq C_n^{\mathrm{flat}}<1$.

Combining the three steps proves \cref{eq:strict-new-comparison}.
\end{proof}

\begin{remark}[Meaning of the upper obstruction]
The last inequality is a limitation of the present \emph{relaxed
constraint}, not a counterexample to P\'{o}lya's conjecture. A proof with
coefficient one must use information not encoded by
\cref{eq:normalized-pointwise-system,eq:coupled-integrated-constraint},
or must quantify the positive slack $\mathcal B_k$ in
\cref{eq:normalized-exact-slack}.
\end{remark}

\subsection{Numerical reference values}

The first coefficient column of \cref{tab:constants} records the
classical Li--Yau value $n/(n+2)$, and the second records the rigorously
defined scalar-envelope coefficient $c_n^{\mathrm{env}}$ from \cite{WangXie}.
The final column is a numerical evaluation of the commutator-coupled quotient at
the trial value of $a$ shown in the preceding column.  The values were
obtained by maximizing first in $u$, numerically integrating the
resulting function of $t$ separately on $(0,1)$ and $(1,a)$, and then
performing the outer one-dimensional search in $a$.  They are rounded
to six decimal places and are included only to indicate the size of the
gain; they are not interval-certified bounds and are not used in any
proof.

\begin{table}[ht]
\centering
\caption{Comparison of the classical Li--Yau individual coefficient,
the uncoupled scalar-envelope coefficient, and numerical trial values
of the commutator-coupled quotient.}\label{tab:constants}
\begin{tabular}{@{}ccccc@{}}
\toprule
$n$ & Li--Yau $\frac{n}{n+2}$ & $c_n^{\mathrm{env}}$
& trial $a$ & coupled $\bigl((a-1)/J_n(a)\bigr)^{2/n}$ \\
\midrule
2 & 0.500000 & 0.538307 & 2.223321 & 0.681188 \\
3 & 0.600000 & 0.622248 & 1.769572 & 0.748347 \\
4 & 0.666667 & 0.681250 & 1.559532 & 0.792643 \\
5 & 0.714286 & 0.724595 & 1.438974 & 0.823897 \\
6 & 0.750000 & 0.757678 & 1.360913 & 0.847074 \\
7 & 0.777778 & 0.783719 & 1.306303 & 0.864922 \\
8 & 0.800000 & 0.804735 & 1.265981 & 0.879079 \\
9 & 0.818182 & 0.822044 & 1.235002 & 0.890574 \\
10 & 0.833333 & 0.836544 & 1.210462 & 0.900090 \\
\bottomrule
\end{tabular}
\end{table}

The table illustrates the strict hierarchy proved in
\cref{prop:constant-comparison}; the last column is diagnostic and is not
used in the proof.  This hierarchy concerns individual eigenvalues.  The
next two sections instead retain the sharp semiclassical coefficient in
the eigenvalue sum and improve its finite-$k$ remainder.

\subsection{A residual vector-field refinement}

There is a natural place where additional geometry could enter. The gradient
space
\[
\mathcal G=\{\nabla u:u\in H_0^1(\Omega)\}
\subset L^2(\Omega;\C^n)
\]
is closed by Poincar\'{e}'s inequality.  Since $\mathbf P_k$ is an
orthogonal projection defined in Section \ref{Section_3.2}, 
$\Ker\mathbf P_k$ is also closed.  Hence
\[
\mathcal H_k=\mathcal G^\perp\cap\Ker\mathbf P_k
\]
is also a closed subspace of $L^2(\Omega;\C^n)$.  We denote by
\[
\Pi_{\mathcal H_k}:L^2(\Omega;\C^n)\longrightarrow \mathcal H_k
%L^2(\Omega;\C^n)
\]
the orthogonal projection onto $\mathcal H_k$, i.e., for every
$Z\in L^2(\Omega;\C^n)$,
\[
\Pi_{\mathcal H_k}Z\in\mathcal H_k,
\qquad
\bigl\langle Z-\Pi_{\mathcal H_k}Z,W\bigr\rangle_{L^2(\Omega;\C^n)}=0
\quad\text{for all }W\in\mathcal H_k.
\]
In particular,
$\Pi_{\mathcal H_k}^2=\Pi_{\mathcal H_k}
=\Pi_{\mathcal H_k}^*$ and
$\Ran\Pi_{\mathcal H_k}=\mathcal H_k$.  Define
\[
\rho_k(\xi)=M_\Omega^{-1}
\|\Pi_{\mathcal H_k}Y_\xi\|_2^2.
\]
Because $\mathcal H_k\subset\Ker\mathbf P_k$, the projection vanishes
on $\Ran\mathbf P_k$, and therefore
\[
\Pi_{\mathcal H_k}Y_\xi
=\Pi_{\mathcal H_k}(I-\mathbf P_k)Y_\xi,
\qquad
0\leq\rho_k(\xi)
\leq M_\Omega^{-1}\|(I-\mathbf P_k)Y_\xi\|_2^2
=1-\frac{F_k(\xi)}{M_\Omega},
\]
where the last equality follows from
\cref{eq:Y-norm,eq:PY-formula}.  Thus the quantity
$1-f-\rho_k$ below is nonnegative; in the normalized variables we use
the abbreviation $\rho_k=\rho_k(\sqrt\Lambda\,\eta)$.  Moreover,
$\Pi_{\mathcal H_k}Y_\xi$ is orthogonal to $\nabla q_\xi$ and, being in
$\Ker\mathbf P_k$, to the componentwise low-mode space
$\Ran\mathbf P_k$.  Removing this component before the second
Cauchy--Schwarz estimate refines
\cref{eq:normalized-pointwise-system} to
\begin{equation}\label{eq:rho-refined-coupling}
\sqrt t\,f
\leq\sqrt{fg}+\sqrt{(1-f-\rho_k)(f-g)}.
\end{equation}
In particular, at $g=0$ one has
\begin{equation}\label{eq:rho-refined-cap}
(1+t)f\leq1-\rho_k.
\end{equation}
This is a genuine improvement whenever $\rho_k$ has a positive lower
bound on a relevant frequency region. No such lower bound follows from
volume alone; accordingly, \cref{eq:rho-refined-coupling} is recorded as
a domain-sensitive refinement and is omitted from the universal
coefficient.

\section{Frequency Concentration and Multifrequency Constraints}
\label{sec:concentration}

The second part replaces the pointwise energy envelope by a collective
constraint on the first $k$ Fourier transforms.  Its central object is a
bounded-domain analogue of the classical time--frequency concentration
operator \cite{Slepian1964}.  We establish its strict spectral structure,
derive first- and second-trace overlap formulas, and then convert them
into a finite-rank restriction for the Fourier projection onto the first
$k$ eigenfunctions.  Throughout \cref{sec:concentration,sec:strict-sum},
$\Omega$ is assumed bounded.

\subsection{The frequency-concentration operator}

Define
\begin{equation}\label{eq:gamma}
\gamma_\Omega(z)=\frac1V\int_\Omega e^{\ii x\cdot z}\dd x.
\end{equation}
It follows directly that $\gamma_\Omega(0)=1$, while the triangle
inequality gives
\[
|\gamma_\Omega(z)|
\leq \frac1V\int_\Omega|e^{\ii x\cdot z}|\dd x
=\frac1V\int_\Omega1\dd x=1.
\]
Let $\one_\Omega$ denote the indicator function of $\Omega$.  Since
$V=|\Omega|<\infty$, we have $\one_\Omega\in L^1(\R^n)\cap L^2(\R^n)$.
With the Fourier convention in \cref{eq:fourier}, the definition in
\cref{eq:gamma} becomes
\[
\gamma_\Omega(z)
=\frac{(2\pi)^{n/2}}{V}\,
\widehat{\one_\Omega}(-z).
\]
Using this identity and $M_\Omega=V/(2\pi)^n$, changing $z$ to $-z$,
and applying Plancherel's theorem
\cite{SteinShakarchi2003}, we obtain
\begin{align*}
M_\Omega\int_{\R^n}|\gamma_\Omega(z)|^2\dd z
&=\frac{V}{(2\pi)^n}\frac{(2\pi)^n}{V^2}
  \int_{\R^n}|\widehat{\one_\Omega}(-z)|^2\dd z\\
&=\frac1V\int_{\R^n}|\widehat{\one_\Omega}(\xi)|^2\dd\xi\\
&=\frac1V\int_{\R^n}|\one_\Omega(x)|^2\dd x
=\frac{|\Omega|}{V}=1,
\end{align*}
which proves
\begin{equation}\label{eq:gamma-parseval}
M_\Omega\int_{\R^n}|\gamma_\Omega(z)|^2\dd z=1.
\end{equation}

For $r>0$, let $H_r=L^2(B_r)$ and $G=L^2(\Omega)$, and define
\begin{equation}\label{eq:Ur}
(U_rf)(x)=(2\pi)^{-n/2}\int_{B_r}e^{\ii x\cdot\xi}f(\xi)\dd\xi.
\end{equation}
Its adjoint is the restricted Fourier transform
\begin{equation}\label{eq:Ur-adjoint}
(U_r^*u)(\xi)=(2\pi)^{-n/2}
\int_\Omega e^{-\ii x\cdot\xi}u(x)\dd x,
\qquad \xi\in B_r.
\end{equation}
Set
\begin{equation}\label{eq:Tr-Ar}
T_r=U_r^*U_r:H_r\to H_r,
\qquad A_r=U_rU_r^*:G\to G.
\end{equation}

\begin{proposition}[Factorization and kernel]\label{prop:concentration-kernel}
The operator $U_r$ is Hilbert--Schmidt, with
\begin{equation}\label{eq:Ur-HS}
\|U_r\|_{\HS}^2=M_\Omega|B_r|=M_\Omega\omega_nr^n.
\end{equation}
Moreover, $T_r$ and $A_r$ are positive, self-adjoint, compact operators,
and
\begin{equation}\label{eq:Tr-kernel}
(T_rf)(\xi)=M_\Omega\int_{B_r}
\gamma_\Omega(\eta-\xi)f(\eta)\dd\eta.
\end{equation}
Finally,
\begin{equation}\label{eq:Tr-contraction}
0\leq T_r\leq I,
\qquad 0\leq A_r\leq I.
\end{equation}
\end{proposition}

\begin{proof}
Write
\[
K_r(x,\xi)=(2\pi)^{-n/2}e^{\ii x\cdot\xi},
\qquad (x,\xi)\in\Omega\times B_r,
\]
for the integral kernel of $U_r$.  Since $\Omega$ and $B_r$ have finite
measure, $K_r\in L^2(\Omega\times B_r)$ and hence $U_r$ is
Hilbert--Schmidt.  More explicitly,
\begin{align*}
\|U_r\|_{\HS}^2
&=\int_\Omega\int_{B_r}|K_r(x,\xi)|^2\dd\xi\dd x\\
&=(2\pi)^{-n}|\Omega|\,|B_r|
=\frac{V}{(2\pi)^n}\omega_nr^n
=M_\Omega\omega_nr^n,
\end{align*}
which proves \cref{eq:Ur-HS}.  We use here the standard kernel
characterization of Hilbert--Schmidt operators; see, for example,
\cite[Chapter~2]{Simon2005}.

For completeness, the adjoint formula can be checked directly.  If
$f\in H_r$ and $u\in G$, Fubini's theorem (or first the calculation on
dense $L^1\cap L^2$ subspaces and then continuity) gives
\begin{align*}
\langle U_rf,u\rangle_G
&=(2\pi)^{-n/2}\int_\Omega\int_{B_r}
e^{\ii x\cdot\eta}f(\eta)\overline{u(x)}\dd\eta\dd x\\
&=\int_{B_r}f(\eta)\overline{
(2\pi)^{-n/2}\int_\Omega e^{-\ii x\cdot\eta}u(x)\dd x}
\dd\eta,
\end{align*}
which is \cref{eq:Ur-adjoint}.  Substituting this formula into
$T_r=U_r^*U_r$ and interchanging the two finite-measure integrals yields
\begin{align*}
(T_rf)(\xi)
&=(2\pi)^{-n}\int_{B_r}\int_\Omega
e^{\ii x\cdot(\eta-\xi)}\dd x\,f(\eta)\dd\eta\\
&=\frac{V}{(2\pi)^n}\int_{B_r}
\left(\frac1V\int_\Omega e^{\ii x\cdot(\eta-\xi)}\dd x\right)
f(\eta)\dd\eta\\
&=M_\Omega\int_{B_r}\gamma_\Omega(\eta-\xi)f(\eta)\dd\eta,
\end{align*}
which proves \cref{eq:Tr-kernel} with all normalization constants displayed.

The identities $T_r=U_r^*U_r$ and $A_r=U_rU_r^*$ show immediately that
both operators are self-adjoint and positive.  Since $U_r$ and $U_r^*$
are Hilbert--Schmidt, these two products are trace class, and in
particular compact.

It remains to verify the operator bounds.  Let $E_r:H_r\to L^2(\R^n)$
denote extension by zero and $R_\Omega:L^2(\R^n)\to L^2(\Omega)$
restriction.  With the unitary Fourier convention used here,
\[
U_r=R_\Omega\mathcal F^{-1}E_r.
\]
Thus $\|U_r\|\leq1$ by Plancherel.  Equivalently, for every $f\in H_r$,
\begin{equation}\label{eq:Tr-quadratic}
\langle T_rf,f\rangle_{H_r}=\|U_rf\|_G^2\leq\|f\|_{H_r}^2.
\end{equation}
Together with positivity this is precisely $0\leq T_r\leq I$.  Likewise,
$U_r^*$ is the Fourier transform of the zero extension of $u\in G$,
followed by restriction to $B_r$, so $\|U_r^*u\|_{H_r}\leq\|u\|_G$.
Therefore
\[
0\leq\langle A_ru,u\rangle_G=\|U_r^*u\|_{H_r}^2
\leq\|u\|_G^2,
\]
which proves $0\leq A_r\leq I$ and completes the proof.
\end{proof}

The strict spectral structure is the essential point.  We use the
following elementary uniqueness statement.

\begin{lemma}[Analytic continuation and uniqueness]
\label{lem:analytic-uniqueness}
Let $g\in L^1(\R^n)$ have bounded support.  Then $\widehat g$ extends
to an entire function on $\C^n$.  Consequently, if $\widehat g$
vanishes on a nonempty open subset of $\R^n$, then $g=0$ almost
everywhere.  The same conclusions hold for the inverse Fourier
transform of a compactly supported $L^2$ function.
\end{lemma}

\begin{proof}
Choose $R>0$ such that $\operatorname{supp}g\subset B_R$.  The function
\[
G(z)=(2\pi)^{-n/2}\int_{B_R}e^{-\ii x\cdot z}g(x)\dd x.
\]
is entire, because the integrand and all of its $z$-derivatives are
locally uniformly dominated by an integrable multiple of $|g|$.  This
is the elementary compact-support case of the Paley--Wiener principle;
see, for example,
\cite[Chapter~4, Section~3]{SteinShakarchi2003Complex}.

If $G$ vanishes on a nonempty open subset of $\R^n$, choose a real open
rectangle contained in that set and apply the one-variable identity
theorem successively in each coordinate.  It follows that
$G\equiv0$ on $\C^n$.  Injectivity of the Fourier transform on
$L^1(\R^n)$ then gives $g=0$ almost everywhere; see
\cite[Chapter~6, Section~2]{SteinShakarchi2003}.

Finally, a compactly supported $L^2$ function belongs to $L^1$ by
Cauchy--Schwarz.  The inverse-transform assertions follow from the same
argument with the sign of the exponential reversed.
\end{proof}

\begin{theorem}[Strict spectral structure]\label{thm:concentration-spectrum}
The nonzero eigenvalues of $T_r$, repeated according to multiplicity,
can be arranged as
\begin{equation}\label{eq:mu-spectrum}
1>\mu_1(r)\geq\mu_2(r)\geq\cdots>0,
\qquad \mu_\ell(r)\longrightarrow0.
\end{equation}
They form the complete spectrum of $T_r$ apart from the accumulation
point zero.  The operator $A_r$ has the same eigenvalues, including
multiplicity, and its corresponding eigenfunctions form an orthonormal
basis of $G$.
\end{theorem}

\begin{proof}
Let $0\neq f\in H_r$, extend $f$ by zero to
$\widetilde f\in L^2(\R^n)$, and put
\[
F=\mathcal F^{-1}\widetilde f.
\]
Because $\widetilde f$ is supported in $\overline{B_r}$, it also belongs
to $L^1(\R^n)$, and $F$ is real analytic (indeed entire after complexification)
by \cref{lem:analytic-uniqueness}.  Moreover,
\[
U_rf=F|_\Omega,
\qquad
\langle T_rf,f\rangle_{H_r}=\|F\|_{L^2(\Omega)}^2,
\qquad
\|F\|_{L^2(\R^n)}=\|f\|_{H_r}.
\]

If the lower endpoint were attained, then $F=0$ almost everywhere on the
nonempty open set $\Omega$.  Since $F$ is continuous, almost-everywhere
vanishing on an open set implies pointwise vanishing there.  Analytic uniqueness would imply
$\widetilde f=0$, contrary to $f\neq0$.  If the upper endpoint were
attained, then
\[
0=\|F\|_{L^2(\R^n)}^2-\|F\|_{L^2(\Omega)}^2
=\|F\|_{L^2(\R^n\setminus\Omega)}^2.
\]
Since $\Omega$ is bounded, $\R^n\setminus\overline\Omega$ is a nonempty
open set; $F$ vanishes there almost everywhere and hence everywhere by
continuity.  Analytic uniqueness again forces $f=0$.  We have therefore proved the strict quadratic-form
bounds
\begin{equation}\label{eq:strict-form}
0<\langle T_rf,f\rangle<\|f\|^2
\qquad(0\neq f\in H_r).
\end{equation}
In particular, $T_r$ is injective: $T_rf=0$ would imply
$\langle T_rf,f\rangle=0$.  It also has no eigenvalue equal to one,
because $T_rf=f$ would imply
$\langle T_rf,f\rangle=\|f\|^2$.  By
\cref{prop:concentration-kernel}, $T_r$ is compact, positive, and
self-adjoint.  The compact self-adjoint spectral theorem
\cite{Teschl2009} therefore shows that every nonzero spectral value is a
positive eigenvalue of finite multiplicity, that the only possible
accumulation point is zero, and that the corresponding eigenspaces span
$(\Ker T_r)^\perp=H_r$.  Since $H_r$ is infinite dimensional whereas an
injective finite-rank operator on $H_r$ cannot exist, there are
infinitely many such eigenvalues.  Listing them with multiplicity gives
\[
\mu_1(r)\geq\mu_2(r)\geq\cdots>0,
\qquad \mu_\ell(r)\to0.
\]
Compactness also implies that a nonzero operator norm is attained as an
eigenvalue; the exclusion of the eigenvalue one therefore yields
$\mu_1(r)=\|T_r\|<1$, proving \cref{eq:mu-spectrum}.

We next compare $T_r$ and $A_r$.  For every $\mu>0$, the map
\[
U_r:\Ker(T_r-\mu I)\longrightarrow\Ker(A_r-\mu I)
\]
is well defined, because
\[
A_r(U_rf)=U_r(U_r^*U_rf)=U_r(T_rf)=\mu U_rf.
\]
It is injective on this eigenspace since
$\|U_rf\|^2=\langle T_rf,f\rangle=\mu\|f\|^2$.  Conversely, if
$A_ru=\mu u$, then
\[
T_r(U_r^*u)=U_r^*(U_rU_r^*u)=\mu U_r^*u,
\]
and $U_r^*u\neq0$ because
$\|U_r^*u\|^2=\langle A_ru,u\rangle=\mu\|u\|^2$.  The inverse of the
first map is $\mu^{-1}U_r^*$, so the two $\mu$-eigenspaces are isomorphic.
Thus $T_r=U_r^*U_r$ and $A_r=U_rU_r^*$ have the same nonzero eigenvalues,
including multiplicities; this is also the standard singular-value
correspondence for a compact operator \cite[Chapter~1]{Simon2005}.

It remains to prove completeness in $G$.  Suppose $U_r^*u=0$.  If
$\widetilde u$ denotes the zero extension of $u$ from $\Omega$ to
$\R^n$, then $\widetilde u$ has bounded support and belongs to $L^1$ by
Cauchy--Schwarz, while
\[
\widehat{\widetilde u}(\xi)=0\qquad\text{for almost every }\xi\in B_r.
\]
The transform is continuous (indeed entire), so it vanishes everywhere
on the open ball $B_r$.  By \cref{lem:analytic-uniqueness}, $u=0$.
Hence $\Ker U_r^*=\{0\}$ and therefore
$\Ker A_r=\Ker U_r^*=\{0\}$.  Applying the compact self-adjoint spectral
theorem to $A_r$ gives
\[
G=\Ker A_r\oplus
\overline{\bigoplus_{\mu>0}\Ker(A_r-\mu I)}
=\overline{\bigoplus_{\mu>0}\Ker(A_r-\mu I)}.
\]
Choosing an orthonormal basis in each finite-dimensional eigenspace
therefore produces the asserted orthonormal eigenbasis of $G$.
\end{proof}

\begin{remark}
The pointwise strict inequality in \cref{eq:strict-form} does not imply
a uniform gap below one: for fixed $\Omega$, $\|T_r\|\to1$ as
$r\to\infty$.  What is used below is the strict inequality for every
fixed finite radius.
\end{remark}

\subsection{Trace identities and overlap asymptotics}\label{sec:traces}

\begin{theorem}[First and second traces]\label{thm:traces}
The operator $T_r$ is positive trace class, $T_r^2$ is trace class, and
\begin{align}
t(r):=\Tr (T_r)
&=M_\Omega|B_r|=M_\Omega\omega_nr^n,
\label{eq:first-trace}\\
Q(r):=\Tr(T_r^2)
&=\sum_{\ell=1}^\infty\mu_\ell(r)^2
=M_\Omega^2\int_{B_r}\int_{B_r}
|\gamma_\Omega(\eta-\xi)|^2\dd\eta\dd\xi.
\label{eq:second-trace}
\end{align}
Moreover,
\begin{equation}\label{eq:trace-gap}
0<Q(r)<t(r).
\end{equation}
\end{theorem}

\begin{proof}
Since $U_r$ is Hilbert--Schmidt, $T_r=U_r^*U_r$ is trace class and
$\Tr(T_r)=\|U_r\|_{\HS}^2$, proving \cref{eq:first-trace}.
The Hilbert--Schmidt norm of the kernel in \cref{eq:Tr-kernel} gives
\cref{eq:second-trace}.  Finally, every $\mu_\ell(r)$ lies in $(0,1)$,
so
$0<\sum_\ell\mu_\ell(r)^2<\sum_\ell\mu_\ell(r)=t(r)$.
\end{proof}

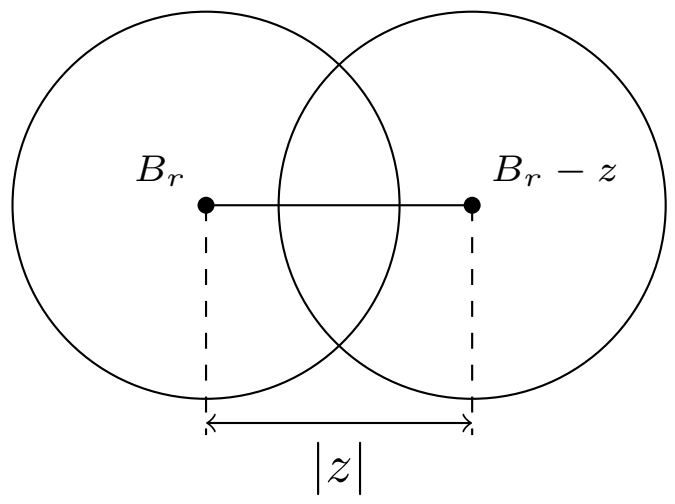
\begin{figure}[ht]
\centering 
\begin{tikzpicture}[scale=0.8]
    % 两个球的圆心
    \coordinate (O1) at (0,0);
    \coordinate (O2) at (2.2,0);
    \def\radius{1.6}
    % 画两个相交圆
    \draw (O1) circle (\radius) node[above left] {\footnotesize $B_r$};
    \draw (O2) circle (\radius) node[above right] {\footnotesize $B_r - z$};

    % 圆心黑点
    \fill (O1) circle (2pt);
    \fill (O2) circle (2pt);
    
    \draw[solid] (O1) -- (O2); 

    % 向量箭头 e_1
    % \draw[->] (O1) -- (1.6,0) node at (1,0.2) {$e_1$};

    % 下方虚线与标注 s
\draw[dashed] (O1) -- (O1 |- {0,-1.9});
\draw[dashed] (O2) -- (O2 |- {0,-1.9});
\draw[<->] (O1 |- {0,-1.8}) -- (O2 |- {0,-1.8}) node[midway, below] {$|z|$};
\end{tikzpicture}
\caption{The overlap of the balls $B_r$ and $B_r-z$ with center
separation $|z|$.}\label{fig:ball-overlap}
\end{figure}

\begin{theorem}[Overlap representation]\label{thm:overlap}
Let
\begin{equation}\label{eq:overlap-function}
a_r(z)=|B_r\cap(B_r-z)|.
\end{equation}
Then
\begin{equation}\label{eq:Q-overlap}
Q(r)=M_\Omega^2\int_{\R^n}
|\gamma_\Omega(z)|^2a_r(z)\dd z,
\end{equation}
and
\begin{equation}\label{eq:Q-asymptotic}
\frac{Q(r)}{t(r)}\longrightarrow 1
\qquad(r\to\infty).
\end{equation}
\end{theorem}
Fig.~\ref{fig:ball-overlap} illustrates the intersection defining $a_r(z)$.

\begin{proof}
Set
\[
\Phi(z)=|\gamma_\Omega(z)|^2.
\]
The integrand in \cref{eq:second-trace} is nonnegative, so Tonelli's
theorem allows all changes in the order of integration below without any
prior integrability assumption.  Extending the integrations by means of
the indicator $\ind_{B_r}$ and setting $z=\eta-\xi$, we obtain
\begin{align*}
\int_{B_r}\int_{B_r}\Phi(\eta-\xi)\dd\eta\dd\xi
&=\int_{\R^n}\int_{\R^n}
  \Phi(\eta-\xi)\ind_{B_r}(\xi)\ind_{B_r}(\eta)
  \dd\eta\dd\xi\\
&=\int_{\R^n}\int_{\R^n}
  \Phi(z)\ind_{B_r}(\xi)\ind_{B_r}(\xi+z)
  \dd\xi\dd z\\
&=\int_{\R^n}\Phi(z)
  \left(\int_{\R^n}
  \ind_{B_r}(\xi)\ind_{B_r-z}(\xi)\dd\xi\right)\dd z\\
&=\int_{\R^n}\Phi(z)|B_r\cap(B_r-z)|\dd z\\
&=\int_{\R^n}|\gamma_\Omega(z)|^2a_r(z)\dd z.
\end{align*}
Multiplication by $M_\Omega^2$ and comparison with
\cref{eq:second-trace} prove \cref{eq:Q-overlap}.  This calculation also
shows directly that $a_r$ is measurable, $0\leq a_r(z)\leq|B_r|$, and
$a_r(z)=0$ when $|z|\geq2r$.

We next justify the pointwise limit needed for the asymptotic formula.
By scaling,
\[
\frac{a_r(z)}{|B_r|}
=\frac{|B_1\cap(B_1-z/r)|}{|B_1|}.
\]
For any vector $h$,
\[
|B_1|-|B_1\cap(B_1-h)|
=\frac12\,|B_1\mathbin{\triangle}(B_1-h)|,
\]
where $\triangle$ denotes symmetric difference.  Translation is
continuous in $L^1(\R^n)$, applied here to $\ind_{B_1}$; hence the last
quantity tends to zero as $h\to0$.  Taking $h=z/r$ gives, for every fixed
$z\in\R^n$,
\[
\frac{a_r(z)}{|B_r|}\longrightarrow1
\qquad(r\to\infty).
\]

Finally, divide \cref{eq:Q-overlap} by
$t(r)=M_\Omega|B_r|$ to get
\begin{equation}\label{eq:normalized-overlap-trace}
\frac{Q(r)}{t(r)}
=M_\Omega\int_{\R^n}|\gamma_\Omega(z)|^2
\frac{a_r(z)}{|B_r|}\dd z. 
\end{equation}
The ratio in the integrand lies between zero and one, while
$M_\Omega|\gamma_\Omega|^2$ is integrable and has integral one by
\cref{eq:gamma-parseval}.  Dominated convergence therefore gives
\[
\lim_{r\to\infty}\frac{Q(r)}{t(r)}
=M_\Omega\int_{\R^n}|\gamma_\Omega(z)|^2\dd z=1,
\]
which is \cref{eq:Q-asymptotic} and the proof is complete. 
\end{proof}

\subsection{Finite-rank projections and the multifrequency constraint}
\label{sec:multifrequency}

Let $P$ be an orthogonal projection of rank $k$ in $G$, and define
\begin{equation}\label{eq:MrP}
\mathcal M_r(P):=\Tr(PA_r).
\end{equation}
If $\{u_j\}_{j=1}^k$ is an orthonormal basis of $\Ran(P)$, then
\begin{equation}\label{eq:MrP-basis}
\mathcal M_r(P)=\sum_{j=1}^k\int_{B_r}
|\widehat u_j(\xi)|^2\dd\xi.
\end{equation}

\begin{theorem}[Ky Fan bound]\label{thm:kyfan}
With
\begin{equation}\label{eq:Bkr}
B_k(r)=\sum_{\ell=1}^k\mu_\ell(r),
\end{equation}
one has
\begin{equation}\label{eq:kyfan}
\mathcal M_r(P)\leq B_k(r),
\qquad
\sup_{\rank P=k}\mathcal M_r(P)=B_k(r).
\end{equation}
Furthermore, we have 
\begin{equation}\label{eq:B-bounds}
B_k(r)<\min\{k,t(r)\},
\qquad B_k(r)\leq\sqrt{kQ(r)}.
\end{equation}
\end{theorem}

\begin{proof}
For brevity write $\mu_\ell=\mu_\ell(r)$.  By
\cref{thm:concentration-spectrum}, $A_r$ has an orthonormal eigenbasis
$\{v_\ell\}_{\ell\geq1}$ with
\[
A_rv_\ell=\mu_\ell v_\ell,
\qquad 1>\mu_1\geq\mu_2\geq\cdots>0.
\]
Since $A_r$ is trace class and $P$ is bounded, $PA_r$ is trace class, so
its trace may be computed in this basis.  Put
\[
a_\ell=\langle Pv_\ell,v_\ell\rangle
=\|Pv_\ell\|^2.
\]
The last equality uses $P=P^*=P^2$.  Thus $0\leq a_\ell\leq1$, and
\begin{align*}
\mathcal M_r(P)
&=\Tr(PA_r)
=\sum_{\ell=1}^\infty\langle PA_rv_\ell,v_\ell\rangle\\
&=\sum_{\ell=1}^\infty
  \mu_\ell\langle Pv_\ell,v_\ell\rangle
=\sum_{\ell=1}^\infty\mu_\ell a_\ell.
\end{align*}

We also record why the coefficients have total mass $k$.  If
$\{u_j\}_{j=1}^k$ is an orthonormal basis of $\Ran(P)$, then
\[
Pv_\ell=\sum_{j=1}^k\langle v_\ell,u_j\rangle u_j.
\]
Parseval's identity in the eigenbasis $\{v_\ell\}$ therefore gives
\begin{align*}
\sum_{\ell=1}^\infty a_\ell
&=\sum_{\ell=1}^\infty\sum_{j=1}^k
  |\langle v_\ell,u_j\rangle|^2
=\sum_{j=1}^k\|u_j\|^2=k.
\end{align*}

The maximizing arrangement places this unit-bounded mass on the largest
$k$ eigenvalues.  Here is a direct verification.  From
$\sum_{\ell>k}a_\ell=k-\sum_{\ell=1}^ka_\ell
=\sum_{\ell=1}^k(1-a_\ell)$ and
$\mu_\ell\leq\mu_k$ for $\ell>k$, we obtain
\begin{align*}
\mathcal M_r(P)
&=\sum_{\ell=1}^k\mu_\ell a_\ell
  +\sum_{\ell>k}\mu_\ell a_\ell\\
&\leq\sum_{\ell=1}^k\mu_\ell a_\ell
  +\mu_k\sum_{\ell=1}^k(1-a_\ell)\\
&=\sum_{\ell=1}^k
  \bigl[\mu_\ell a_\ell+\mu_k(1-a_\ell)\bigr]\\
&\leq\sum_{\ell=1}^k\mu_\ell=B_k(r).
\end{align*}
Let $P_k^{\mathrm{sp}}$ be the orthogonal projection onto
$\operatorname{span}\{v_1,\ldots,v_k\}$.  For this projection,
$a_\ell=1$ for $\ell\leq k$ and $a_\ell=0$ for $\ell>k$, so
$\mathcal M_r(P_k^{\mathrm{sp}})=B_k(r)$.  This proves both assertions in
\cref{eq:kyfan}, including attainment of the supremum.

It remains to prove \cref{eq:B-bounds}.  Since every $\mu_\ell<1$,
\[
B_k(r)=\sum_{\ell=1}^k\mu_\ell<k.
\]
Moreover, \cref{thm:concentration-spectrum} gives infinitely many
strictly positive eigenvalues, and \cref{eq:first-trace} gives
$t(r)=\sum_{\ell=1}^\infty\mu_\ell$.  Hence
\[
t(r)-B_k(r)=\sum_{\ell>k}\mu_\ell>0,
\]
which proves $B_k(r)<t(r)$.  Finally, Cauchy--Schwarz and
\cref{eq:second-trace} yield
\[
B_k(r)^2
\leq k\sum_{\ell=1}^k\mu_\ell^2
\leq k\sum_{\ell=1}^\infty\mu_\ell^2
=kQ(r).
\]
Taking square roots completes the proof.
\end{proof}

Let $P_k$ be the projection onto
$\operatorname{span}\{\phi_1,\ldots,\phi_k\}$.  By
\cref{eq:density,eq:MrP-basis,eq:kyfan,eq:B-bounds}, we have 
\begin{equation}\label{eq:inside-mass}
\int_{B_r}F_k(\xi)\dd\xi
\leq B_k(r)\leq\sqrt{kQ(r)}.
\end{equation}
This collective constraint is absent from the pointwise envelopes in
\cref{eq:l2-bessel,eq:projection-envelope}.

To convert the resulting mass leakage into a moment estimate, we need a
two-region version of the bathtub principle.

\begin{lemma}[Minimal moment at prescribed outside mass]
\label{lem:leakage-bathtub}
Let $M_0,R>0$, set $k=M_0\omega_nR^n$, and let
$f:\R^n\to[0,M_0]$ satisfy $\int_{\R^n}f=k$.  If
\[
m=\int_{\R^n\setminus B_R}f(\xi)\dd\xi,
\qquad \tau=\frac{m}{k}\in[0,1],
\]
then, with $q=1+2/n$, we have
\begin{equation}\label{eq:leakage-bathtub}
\int_{\R^n}|\xi|^2f(\xi)\dd\xi
\geq\frac{n}{n+2}kR^2
\bigl[(1-\tau)^q+(1+\tau)^q-1\bigr].
\end{equation}
The right-hand side is optimal under these constraints.
\end{lemma}

\begin{proof}
If $\int_{\R^n}|\xi|^2f(\xi)\dd\xi=+\infty$, the estimate is immediate,
so we may assume that the second moment is finite.  Since $f\geq0$ and
its total mass is $k$, one has $0\leq m\leq k$ and hence
$0\leq\tau\leq1$, as stated.

The prescribed masses in the two regions are
\begin{align}
\int_{B_R}f(\xi)\dd\xi
&=k-m=k(1-\tau),
\label{eq:leakage-inner-mass}\\
\int_{\R^n\setminus B_R}f(\xi)\dd\xi
&=m=k\tau.
\label{eq:leakage-outer-mass}
\end{align}
Set
\[
\rho=R(1-\tau)^{1/n},
\qquad \sigma=R(1+\tau)^{1/n}.
\]
Thus $0\leq\rho\leq R\leq\sigma\leq2^{1/n}R$.  The two radii are chosen
so that a density saturated at the cap $M_0$ fills exactly the required
inner and outer masses.
Since $k=M_0\omega_nR^n$, these definitions give
\begin{align*}
M_0|B_\rho|
&=M_0\omega_nR^n(1-\tau)=k(1-\tau),\\
M_0|B_\sigma\setminus B_R|
&=M_0\omega_n(\sigma^n-R^n)=k\tau.
\end{align*}
Thus $M_0\ind_{B_\rho}$ has exactly the same mass as $f$ in $B_R$,
and $M_0\ind_{B_\sigma\setminus B_R}$ has exactly the same mass as
$f$ in $\R^n\setminus B_R$.

The least expensive inner configuration is
$M_0\ind_{B_\rho}$, while the least expensive outer configuration is
$M_0\ind_{B_\sigma\setminus B_R}$.  To see the first assertion,
observe that on $B_\rho$ both factors in
\[
(|\xi|^2-\rho^2)
\bigl(f-M_0\ind_{B_\rho}\bigr)\geq0,
\]
are nonpositive, whereas on $B_R\setminus B_\rho$ they are both
nonnegative.  Hence the displayed inequality holds throughout $B_R$.
To display the cancellation rather than hide it in the rearrangement
argument, write
\begin{align*}
&\int_{B_R}|\xi|^2
  \bigl(f-M_0\ind_{B_\rho}\bigr)\dd\xi\\
&\quad=\int_{B_R}(|\xi|^2-\rho^2)
  \bigl(f-M_0\ind_{B_\rho}\bigr)\dd\xi
  +\rho^2\left(\int_{B_R}f\dd\xi-M_0|B_\rho|\right).
\end{align*}
The first term is nonnegative by the preceding pointwise sign argument,
whereas the term in parentheses vanishes by
\cref{eq:leakage-inner-mass}.  It follows that
\begin{equation}\label{eq:leakage-inner-moment}
\int_{B_R}|\xi|^2f(\xi)\dd\xi
\geq M_0\int_{B_\rho}|\xi|^2\dd\xi.
\end{equation}

Similarly, on $B_\sigma\setminus B_R$ both factors in
\[
(|\xi|^2-\sigma^2)
\bigl(f-M_0\ind_{B_\sigma\setminus B_R}\bigr)\geq0,
\]
are nonpositive, while on $\R^n\setminus B_\sigma$ both are
nonnegative.  The corresponding expanded identity is
\begin{align*}
&\int_{\R^n\setminus B_R}|\xi|^2
  \bigl(f-M_0\ind_{B_\sigma\setminus B_R}\bigr)\dd\xi\\
&\quad=\int_{\R^n\setminus B_R}(|\xi|^2-\sigma^2)
  \bigl(f-M_0\ind_{B_\sigma\setminus B_R}\bigr)\dd\xi\\
&\qquad+\sigma^2\left(
  \int_{\R^n\setminus B_R}f\dd\xi
  -M_0|B_\sigma\setminus B_R|\right).
\end{align*}
The first term is nonnegative by the sign argument and the second term
vanishes by \cref{eq:leakage-outer-mass}.  Therefore
\begin{equation}\label{eq:leakage-outer-moment}
\int_{\R^n\setminus B_R}|\xi|^2f(\xi)\dd\xi
\geq M_0\int_{B_\sigma\setminus B_R}|\xi|^2\dd\xi.
\end{equation}

For completeness, polar coordinates give, for every $a\geq0$,
\begin{equation}\label{eq:ball-second-moment}
\int_{B_a}|\xi|^2\dd\xi
=|\mathbb S^{n-1}|\int_0^a r^{n+1}\dd r
=n\omega_n\frac{a^{n+2}}{n+2}.
\end{equation}
Adding \cref{eq:leakage-inner-moment,eq:leakage-outer-moment} and then
using \cref{eq:ball-second-moment}, we obtain
\begin{align*}
\int_{\R^n}|\xi|^2f(\xi)\dd\xi
&\geq\frac{M_0n\omega_n}{n+2}
\bigl(\rho^{n+2}+\sigma^{n+2}-R^{n+2}\bigr)\\
&=\frac{M_0n\omega_nR^{n+2}}{n+2}
\left[(1-\tau)^{(n+2)/n}
+(1+\tau)^{(n+2)/n}-1\right]\\
&=\frac{n}{n+2}kR^2
\left[(1-\tau)^q+(1+\tau)^q-1\right],
\end{align*}
where the last line uses $M_0\omega_nR^n=k$ and
$q=(n+2)/n=1+2/n$.  This is \cref{eq:leakage-bathtub}.

Finally,
\[
f_*(\xi)=M_0\ind_{B_\rho}(\xi)
+M_0\ind_{B_\sigma\setminus B_R}(\xi)
\]
takes values in $[0,M_0]$, has total mass $k$ and outside mass $m$,
and makes both regional comparisons equalities.  It therefore attains
the right-hand side of \cref{eq:leakage-bathtub}, proving optimality
rather than merely a lower estimate.  The construction also covers the
endpoints without a limiting argument: if $\tau=0$, then
$\rho=\sigma=R$ and the outer annulus is empty; if $\tau=1$, then
$\rho=0$ and all mass is placed in
$B_{2^{1/n}R}\setminus B_R$.
\end{proof}

\section{Strict Li--Yau Sum Improvements}\label{sec:strict-sum}
We now convert the multifrequency leakage into an eigenvalue-sum
estimate.  We derive the exact spectral and squared-trace multipliers,
prove their quantitative behavior in $k$, and record the resulting
explicit linear correction.

Define the Li--Yau filling radius by
\begin{equation}\label{eq:Rk}
R_k=\left(\frac{k}{M_\Omega\omega_n}\right)^{1/n},
\qquad t(R_k)=k,
\end{equation}
and set
\begin{equation}\label{eq:beta-theta-def}
\beta_k=\frac{B_k(R_k)}{k},
\qquad
\theta_k=\left(\frac{Q(R_k)}{k}\right)^{1/2}.
\end{equation}
By \cref{eq:trace-gap,eq:B-bounds},
\begin{equation}\label{eq:beta-theta-range}
0<\beta_k\leq\theta_k<1.
\end{equation}
At $r=R_k$, \cref{eq:inside-mass} therefore implies
\begin{equation}\label{eq:outside-lower}
\frac1k\int_{\R^n\setminus B_{R_k}}F_k(\xi)\dd\xi
\geq1-\beta_k\geq1-\theta_k>0.
\end{equation}
\begin{theorem}[Exact and computable strict improvements]
\label{thm:strict-sum}
Let $\Omega\subset\R^n$, $n\geq2$, be a nonempty bounded open set.  For
every $k\geq1$, put $q=1+2/n$ and define
\begin{align}
\mathcal C_k^{\mathrm{sp}}(\Omega)
&=\beta_k^q+(2-\beta_k)^q-1,
\label{eq:Csp}\\
\mathcal C_k^{(2)}(\Omega)
&=\theta_k^q+(2-\theta_k)^q-1.
\label{eq:C2}
\end{align}
Then
\begin{equation}\label{eq:strict-sum-chain}
\sum_{j=1}^k\lambda_j
\geq\mathcal C_k^{\mathrm{sp}}(\Omega)L_k^{\mathrm{LY}}
\geq\mathcal C_k^{(2)}(\Omega)L_k^{\mathrm{LY}}
>L_k^{\mathrm{LY}},
\end{equation}
where, as in the introduction,
\begin{equation}\label{eq:LY-symbol}
L_k^{\mathrm{LY}}
=\frac{n}{n+2}C_nV^{-2/n}k^{1+2/n}.
\end{equation}
For fixed $\Omega$,
\begin{equation}\label{eq:C-limits}
\mathcal C_k^{\mathrm{sp}}(\Omega)\longrightarrow1,
\qquad
\mathcal C_k^{(2)}(\Omega)\longrightarrow1
\qquad(k\to\infty).
\end{equation}
\end{theorem}

\begin{proof}
Apply \cref{lem:leakage-bathtub} to $f=F_k$, $M_0=M_\Omega$,
$R=R_k$.  The function
\[
H_n(\tau)=(1-\tau)^q+(1+\tau)^q-1
\]
is strictly increasing on $(0,1)$, because
\[
H_n'(\tau)=q\bigl[(1+\tau)^{q-1}-(1-\tau)^{q-1}\bigr]>0.
\]
In view of
\cref{eq:outside-lower}, we have $\tau\geq 1-\beta_k\geq 1-\theta_k$ and hence
\[
\sum_{j=1}^k\lambda_j
\geq H_n(1-\beta_k)L_k^{\mathrm{LY}}
\geq H_n(1-\theta_k)L_k^{\mathrm{LY}},
\]
which gives the first two inequalities in
\cref{eq:strict-sum-chain}.  Since $0<\theta_k<1$ and $q>1$, strict
convexity yields
\[
\theta_k^q+(2-\theta_k)^q>2.
\]
Thus $\mathcal C_k^{(2)}(\Omega)>1$, proving that the Li--Yau bound is
strictly improved for every finite $k$.

It remains to verify the limits.  Since $R_k\to\infty$ as $k\to\infty$,
\cref{eq:Q-asymptotic} gives $\theta_k^2=Q(R_k)/t(R_k)\to1$.
To treat $\beta_k$, we use the fractional version of the weight-filling
argument from the proof of \cref{thm:kyfan}.  If
$\mu_1\geq\mu_2\geq\cdots\geq0$ and a sequence
$(b_\ell)$ satisfies
\begin{align}
0\leq b_\ell\leq1,
\qquad \sum_{\ell=1}^\infty b_\ell=k,\label{equality_6_1}
\end{align}
then the total available weight $k$ is maximally paired with the $k$
largest values of the decreasing sequence $(\mu_\ell)$.  We give the
full calculation.  First define the amount of weight placed beyond the
first $k$ indices by
\begin{align}
\delta:=\sum_{\ell>k}b_\ell.\label{equality_6_2}
\end{align}
The mass constraint \eqref{equality_6_1} on $(b_\ell)$ implies
\begin{align*}
\delta
&=k-\sum_{\ell=1}^k b_\ell
=\sum_{\ell=1}^k(1-b_\ell).
\end{align*}
Thus $\delta$ is also exactly the unfilled capacity among the first $k$
indices.  Since \eqref{equality_6_2} and $\mu_\ell\leq\mu_k$ for $\ell>k$, the tail contribution
is bounded by
\begin{align}
\sum_{\ell>k}\mu_\ell b_\ell
\leq\mu_k\sum_{\ell>k}b_\ell
=\mu_k\delta.\label{equality_6_3}
\end{align}
On the other hand, $\mu_\ell\geq\mu_k$ for $1\leq\ell\leq k$, so the
unfilled part of the first $k$ positions satisfies
\begin{align}
\sum_{\ell=1}^k\mu_\ell(1-b_\ell)
\geq\mu_k\sum_{\ell=1}^k(1-b_\ell)
=\mu_k\delta.\label{equality_6_4}
\end{align}
Combining the last two estimates \eqref{equality_6_3} and \eqref{equality_6_4} gives
\begin{align}
\sum_{\ell>k}\mu_\ell b_\ell
\leq\sum_{\ell=1}^k\mu_\ell(1-b_\ell).\label{equality_6_5}
\end{align}
Equivalently, from \eqref{equality_6_3}, \eqref{equality_6_4} and \eqref{equality_6_5}, 
subtracting the weighted sum from the completely filled
first $k$ positions yields the explicitly nonnegative deficit
\begin{align*}
\sum_{\ell=1}^k\mu_\ell
-\sum_{\ell=1}^\infty\mu_\ell b_\ell
&=\sum_{\ell=1}^k\mu_\ell(1-b_\ell)
  -\sum_{\ell>k}\mu_\ell b_\ell\\
&\geq \mu_k\delta-\mu_k\delta=0.
\end{align*}
This proves
\begin{equation}\label{eq:fractional-kyfan}
\sum_{\ell=1}^\infty\mu_\ell b_\ell
\leq\sum_{\ell=1}^k\mu_\ell.
\end{equation}
All sums in this argument are finite or absolutely convergent: indeed,
$0\leq\sum_\ell\mu_\ell b_\ell\leq\mu_1\sum_\ell b_\ell
=k\mu_1$.  This also justifies splitting the series at the index $k$.

At $r=R_k$, the first-trace identity and \cref{eq:Rk} give
\[
0<\mu_\ell(R_k)<1,
\qquad
\sum_{\ell=1}^\infty\mu_\ell(R_k)
=t(R_k)=k.
\]
Thus $b_\ell:=\mu_\ell(R_k)$ is an admissible fractional weight sequence
in \cref{eq:fractional-kyfan}: the pointwise constraint
$0<b_\ell<1$ follows from \cref{thm:concentration-spectrum}, and the
preceding display verifies its total mass.  This is only a choice of test
weights; it does not assert that
$b_\ell=\|Pv_\ell\|^2$ for some rank-$k$ projection.  With this choice,
the left- and right-hand sides of \cref{eq:fractional-kyfan} become,
respectively,
\[
\sum_{\ell=1}^\infty\mu_\ell(R_k)b_\ell
=\sum_{\ell=1}^\infty\mu_\ell(R_k)^2=Q(R_k)
\]
and
\[
\sum_{\ell=1}^k\mu_\ell(R_k)=B_k(R_k).
\]
Therefore \cref{eq:fractional-kyfan} gives
$Q(R_k)\leq B_k(R_k)$.  Dividing by $k=t(R_k)$ and using
\cref{eq:beta-theta-def} yields
\[
\theta_k^2=\frac{Q(R_k)}k
\leq\frac{B_k(R_k)}k=\beta_k.
\]
Together with the already proved bound $\beta_k\leq\theta_k$, this gives
$\theta_k^2\leq\beta_k\leq\theta_k$.  Since $\theta_k\to1$, the squeeze
theorem now shows $\beta_k\to1$ as well.
\end{proof}

The squared-trace version is explicit in terms of the Fourier transform
of the domain.  Indeed,
\cref{eq:normalized-overlap-trace,eq:beta-theta-def} gives
\begin{equation}\label{eq:theta-direct}
\theta_k^2
=M_\Omega\int_{\R^n}|\gamma_\Omega(z)|^2
\frac{|B_{R_k}\cap(B_{R_k}-z)|}{|B_{R_k}|}\dd z.
\end{equation}
The overlap vanishes for $|z|\geq2R_k$, so only a bounded frequency
region occurs in this formula.

The overlap representation contains additional information about the
dependence on $k$.  The next proposition makes this dependence explicit.

\begin{proposition}[Monotonicity and a quantitative excess]
\label{prop:C2-monotone-gap}
Let $\Omega\subset\R^n$, $n\geq2$, be a nonempty bounded open set.  
Then, we have 
\begin{equation}\label{eq:theta-monotone}
0<\theta_k<\theta_{k+1}<1,
\qquad
\mathcal C_{k+1}^{(2)}(\Omega)
<\mathcal C_k^{(2)}(\Omega).
\end{equation}
Moreover, for every $k\geq1$,
\begin{equation}\label{eq:C2-universal-gap}
\mathcal C_k^{(2)}(\Omega)-1
\geq\kappa_n k^{-2/n},
\end{equation}
where
\begin{equation}\label{eq:kappa-n}
\kappa_n=\frac{n+2}{2(n+1)^2}
\left(\frac{\omega_{n-1}}{\omega_n}\right)^2
\left(\frac34\right)^{n-1}.
\end{equation}
In particular, $\mathcal C_k^{(2)}(\Omega)$ has the following asymptotic property 
\begin{equation}\label{eq:C2-downarrow}
\mathcal C_k^{(2)}(\Omega)\downarrow1.
\end{equation}
\end{proposition}

\begin{proof}
% Introduce the normalized overlap profile (see Fig. \ref{Figure_twoballs})
Rescaling the configuration in Fig.~\ref{fig:ball-overlap} by $1/r$
gives two unit balls with center separation $s=|z|/r$.
By rotational invariance, we may take them to be $B_1$ and $B_1-se_1$.
The normalized overlap volume is then given by
\begin{equation}\label{eq:Phi-overlap}
\Phi_n(s)=\frac{|B_1\cap(B_1-se_1)|}{\omega_n}
=\frac{2\omega_{n-1}}{\omega_n}
\int_{s/2}^{1}(1-u^2)^{(n-1)/2}\dd u,
\qquad 0\leq s\leq2,
\end{equation}
and set $\Phi_n(s)=0$ for $s\geq2$.  Thus
\[
\frac{|B_r\cap(B_r-z)|}{|B_r|}
=\Phi_n\left(\frac{|z|}{r}\right).
\]
% \begin{figure}[ht]
% \centering 
% \begin{tikzpicture}[scale=0.8]
%     % 两个球的圆心
%     \coordinate (O1) at (0,0);
%     \coordinate (O2) at (2.2,0);
%     \def\radius{1.6}
%     % 画两个相交圆
%     \draw (O1) circle (\radius) node[above left] {\footnotesize $B_1$};
%     \draw (O2) circle (\radius) node[above right] {\footnotesize $B_1 - s e_1$};

%     % 圆心黑点
%     \fill (O1) circle (2pt);
%     \fill (O2) circle (2pt);
    
%     \draw[solid] (O1) -- (O2); 

%     % 向量箭头 e_1
%     \draw[->] (O1) -- (1.6,0) node at (1,0.2) {$e_1$};

%     % 下方虚线与标注 s
% \draw[dashed] (O1) -- (O1 |- {0,-1.9});
% \draw[dashed] (O2) -- (O2 |- {0,-1.9});
% \draw[<->] (O1 |- {0,-1.8}) -- (O2 |- {0,-1.8}) node[midway, below] {$s$};
% \end{tikzpicture}
% \caption{Schematic of two intersecting unit balls}\label{Figure_twoballs}
% \end{figure}

The function $\Phi_n$ is strictly decreasing on $(0,2)$.
Consequently, the function
\begin{equation}\label{eq:q-r-def}
q(r):=\frac{Q(r)}{t(r)}
=M_\Omega\int_{\R^n}|\gamma_\Omega(z)|^2
\Phi_n\left(\frac{|z|}{r}\right)\dd z
\end{equation}
is nondecreasing in $r$.  The increase is strict: since
$\gamma_\Omega(0)=1$ and $\gamma_\Omega$ is continuous, its modulus is
bounded away from zero in a neighborhood of the origin, where the
overlap profile increases strictly when $r$ increases.  Since
$R_{k+1}>R_k$ and $\theta_k^2=q(R_k)$, this proves the first inequality
in \cref{eq:theta-monotone}.

Now put
\begin{equation}\label{eq:rho-domain-fourier}
\rho_\Omega(z)=M_\Omega|\gamma_\Omega(z)|^2.
\end{equation}
By \cref{eq:gamma-parseval} and $|\gamma_\Omega|\leq1$,
\begin{equation}\label{eq:rho-cap}
0\leq\rho_\Omega\leq M_\Omega,
\qquad
\int_{\R^n}\rho_\Omega(z)\dd z=1.
\end{equation}
The radial cost
\[
z\longmapsto1-\Phi_n\left(\frac{|z|}{R_k}\right)
\]
is nondecreasing in $|z|$.  The bathtub principle applied to
\cref{eq:rho-cap} therefore shows that its integral is minimized by
$M_\Omega\ind_{B_{R_1}}$, where
$M_\Omega\omega_nR_1^n=1$.  Since $R_k=k^{1/n}R_1$, we obtain the
domain-independent estimate
\begin{align}
1-\theta_k^2
&\geq d_{n,k}
:=n\int_0^1
\left[1-\Phi_n\left(\frac{s}{k^{1/n}}\right)\right]
s^{n-1}\dd s.
\label{eq:exact-deficit-lower}
\end{align}
For $0\leq t\leq1$, \cref{eq:Phi-overlap} gives
\begin{align}
1-\Phi_n(t)
&=\frac{2\omega_{n-1}}{\omega_n}
\int_0^{t/2}(1-u^2)^{(n-1)/2}\dd u
\notag\\
&\geq
\frac{\omega_{n-1}}{\omega_n}
\left(\frac34\right)^{(n-1)/2}t.
\label{eq:overlap-linear-loss}
\end{align}
Since $s/k^{1/n}\leq1$, inserting this estimate into
\cref{eq:exact-deficit-lower} yields
\begin{equation}\label{eq:theta-deficit-rate}
1-\theta_k^2
\geq A_n k^{-1/n},
\qquad
A_n=\frac{n}{n+1}
\frac{\omega_{n-1}}{\omega_n}
\left(\frac34\right)^{(n-1)/2}.
\end{equation}

It remains to pass from the trace deficit to the multiplier.  Let
$\delta=1-\theta_k$ and recall $q=1+2/n\in(1,2]$.  Convexity gives
\begin{equation}\label{eq:symmetric-convexity}
(1-\delta)^q+(1+\delta)^q-2
\geq q(q-1)\delta^2,
\qquad 0\leq\delta\leq1.
\end{equation}
Indeed, the second derivative of the left-hand side is
\[
q(q-1)\bigl[(1-\delta)^{q-2}+(1+\delta)^{q-2}\bigr]
\geq2q(q-1).
\]
For $0\leq\delta<1$, the last inequality follows from AM--GM because
$(1-\delta^2)^{q-2}\geq1$ when $q\leq2$; the endpoint follows by
continuity (or by a one-sided limit when $q<2$).
Since
\[
1-\theta_k=\frac{1-\theta_k^2}{1+\theta_k}
\geq\frac12(1-\theta_k^2),
\]
\cref{eq:theta-deficit-rate,eq:symmetric-convexity} imply
\[
\mathcal C_k^{(2)}(\Omega)-1
\geq\frac{q(q-1)}4A_n^2k^{-2/n}
=\kappa_n k^{-2/n}.
\]
Finally, the function
\[
G(s)=s^q+(2-s)^q-1
\]
satisfies $G'(s)<0$ on $(0,1)$.  The strict increase of $\theta_k$
therefore proves the second assertion in \cref{eq:theta-monotone}; its
limit, already established in \cref{eq:C-limits}, proves
\cref{eq:C2-downarrow}.
\end{proof}

Combining the quantitative excess with \cref{eq:strict-sum-chain}
gives a completely explicit correction to the Berezin--Li--Yau sum.

\begin{corollary}[A universal explicit sum correction]
\label{cor:explicit-sum-correction}
Under the assumptions of \cref{prop:C2-monotone-gap},
\begin{align}
\sum_{j=1}^k\lambda_j
&\geq\left(1+\kappa_n k^{-2/n}\right)L_k^{\mathrm{LY}}
\notag\\
&=L_k^{\mathrm{LY}}
+\frac{n}{n+2}\kappa_nC_nV^{-2/n}k.
\label{eq:explicit-linear-correction}
\end{align}
\end{corollary}

\begin{remark}[The exact spectral coefficient]
\label{rem:Csp-monotonicity}
The same conclusion is not currently available for
$\mathcal C_k^{\mathrm{sp}}$.  Since $G$ is strictly decreasing,
monotonic decrease of $\mathcal C_k^{\mathrm{sp}}$ is equivalent to
$\beta_{k+1}\geq\beta_k$, or
\begin{equation}\label{eq:missing-cross-k}
B_{k+1}(R_{k+1})
\geq\frac{k+1}{k}B_k(R_k).
\end{equation}
For fixed $k$, $B_k(r)$ increases with $r$, whereas for fixed $r$ the
averages $B_k(r)/k$ decrease with $k$.  These two effects compete along
$r=R_k$, and the Ky Fan argument alone does not imply
\cref{eq:missing-cross-k}.  What is proved without an additional
cross-index estimate is
\[
\mathcal C_k^{\mathrm{sp}}(\Omega)>1,
\qquad
\mathcal C_k^{\mathrm{sp}}(\Omega)\longrightarrow1.
\]
Thus no monotonicity claim for the exact spectral coefficient is made.
\end{remark}

\begin{corollary}[The two-dimensional quantitative sum bound]
\label{cor:sum-2d}
If $n=2$, then
\begin{equation}\label{eq:sum-2d}
\sum_{j=1}^k\lambda_j
\geq\bigl[1+2(1-\theta_k)^2\bigr]\frac{2\pi k^2}{V}
\geq\left(1+\frac{2}{3\pi^2k}\right)
\frac{2\pi k^2}{V}
=\frac{2\pi k^2}{V}+\frac{4k}{3\pi V}.
\end{equation}
\end{corollary}

\begin{remark}[Domain assumptions for the quantitative gap]
The concentration-operator results are stated for bounded open sets,
as in \cref{sec:concentration}.  Once the leakage inequality
has been established, however, the proof of \cref{eq:C2-universal-gap}
uses only $0<V<\infty$, Parseval's identity, $|\gamma_\Omega|\leq1$,
and the Euclidean ball-overlap profile.  Thus no boundary smoothness,
Lipschitz regularity, convexity, connectedness, or finite-perimeter
assumption enters the quantitative gap estimate itself.
\end{remark}

\begin{remark}[Two distinct improvements]
The individual coefficient $C_n^{\mathrm{new}}$ is universal and remains
strictly above the uncoupled envelope coefficient $c_n^{\mathrm{env}}$.
It therefore remains strictly above $n/(n+2)$ as $k\to\infty$. By contrast,
$\mathcal C_k^{(2)}(\Omega)$ improves the sharp Berezin--Li--Yau sum coefficient
for each finite $k$, decreases strictly to one, and has the universal
excess quantified in \cref{eq:C2-universal-gap}.  The two statements
therefore improve different losses in the classical argument and should
not be numerically identified.
\end{remark}

\begin{proposition}[Rigid-motion and scaling behavior]
\label{prop:scaling}
Under the assumptions of \cref{thm:strict-sum}, the coefficients
$\mathcal C_k^{\mathrm{sp}}$ and $\mathcal C_k^{(2)}$ are invariant
under translations and uniform scalings of the domain.
\end{proposition}

\begin{proof}
If $\Omega'=a\Omega+x_0$, $a>0$, then
$\gamma_{\Omega'}(z)=e^{\ii x_0\cdot z}\gamma_\Omega(az)$ and
$R_k(\Omega')=a^{-1}R_k(\Omega)$.  A change of variables in
\cref{eq:Q-overlap} proves the assertion for $\mathcal C_k^{(2)}$; the
concentration operators are unitarily equivalent after rescaling, which
proves it for $\mathcal C_k^{\mathrm{sp}}$.
\end{proof}

\section{Conclusion and Outlook}\label{sec:questions}

This paper establishes two complementary refinements of the
Berezin--Li--Yau Fourier method.  For individual eigenvalues, the exact
commutator identity in \cref{prop:exact-commutator-slack} and the coupled
scalar--vector projection inequality in \cref{prop:coupled-projection}
lead to \cref{thm:commutator-coupled-main}.  The comparison in
\cref{prop:constant-comparison,tab:constants} makes explicit the two
successive gains: the scalar energy envelope first improves the
Li--Yau coefficient, and the coupled constraint then improves that
envelope by retaining compatibility between the scalar Fourier mass and
its projected gradient.

For eigenvalue sums on bounded open sets, the trace identities and Ky
Fan estimate in \cref{thm:traces,thm:kyfan} impose a collective
frequency constraint on the first $k$ eigenfunctions.  Together with
the leakage estimate, these yield the strict finite-index multiplier in
\cref{thm:strict-sum}; \cref{prop:C2-monotone-gap} quantifies its
universal excess, and \cref{cor:explicit-sum-correction} converts that
excess into a volume-only linear correction.  This part requires no
boundary smoothness, convexity, connectedness, or moment-of-inertia
information.

The limitations are equally clear.  The individual coefficient obtained
from the relaxed coupled constraints is strictly below one, whereas the
sum multiplier decreases to one; neither result therefore proves
P\'{o}lya's conjecture.  Three structures retained by the analysis but
not yet exploited quantitatively are the nonnegative commutator slack
$\mathcal B_k$, the residual vector field $\rho_k$ in
\cref{eq:rho-refined-coupling}, and the concentration spectrum beyond
its first two traces.  Lower bounds for the first two quantities, or
higher-trace and nonradial estimates for the third, are natural next
steps.  To alter the P\'{o}lya coefficient, such refinements would need
signed information at the boundary scale or a mechanism that survives
the passage from averaged spectral bounds to individual eigenvalues.

\backmatter

\section*{Statements and Declarations}

\noindent\textbf{Funding.}
This work was supported by the National Key Research and Development
Program of China (2025YFA1016600, 2025YFA1016601),
the National Natural Science Foundation of China (12331015),
and the National Center for Mathematics and Interdisciplinary Science, CAS.

\noindent\textbf{Competing interests.}
The authors declare that there is no conflict of interest.

\noindent\textbf{Data availability.}
This article is theoretical and does not use external datasets.

\end{document}